\documentclass[10pt]{amsart}
\usepackage[all]{xy}
\usepackage{amssymb}

\newcommand{\cL}{{\mathcal{L}}}

\newcommand{\Spe}{{\rm Spec}}

\newcommand{\Dgl}{{D_{\rm gl}}}

\newcommand{\A}{{ \rm Aut }}

\renewcommand{\O}{{\mathcal{O}}}
\newcommand{\Z}{{\mathbb{Z}}}

\renewcommand{\mod}{{\;\rm mod}}

\newtheorem{pro}{Proposition}[section]
\newtheorem{lemma}[pro]{Lemma}

\newtheorem{cor}[pro]{Corollary}
\newtheorem{theorem}[pro]{Theorem}
\newtheorem{crit}[pro]{Criterion}
\newtheorem{orism}[pro]{Definition}

\begin{document}
\bibliographystyle{amsplain}

\title{On the Krull dimension of the deformation ring of curves with automorphisms}

\author{A. Kontogeorgis}

\begin{abstract} 
We reduce the study of the   Krull dimension $d$ of the deformation ring 
of the functor of deformations of curves with automorphisms
to the study of the tangent space of  the deformation 
functor of a class of matrix representations of the $p$-part 
of the decomposition groups at wild ramified points, 
and we give a method in order to compute $d$. 
\end{abstract}

\email{kontogar@aegean.gr}
\address{
Department of Mathematics, University of the \AE gean, 83200 Karlovassi, Samos,
Greece\\ { \texttt{\upshape http://eloris.samos.aegean.gr}}
}
\date{\today}

\maketitle
\section{Introduction}
Let $X$ be a non-singular  projective curve defined over the field $k$, and let $G$ be a 
 fixed subgroup of the automorphism group of $X$. 
We will denote by  $(X,G)$  the couple of the curve $X$ together with the  
group $G$. 

A deformation of the couple $(X,G)$ over a local ring $A$ is a proper, smooth family
of curves 
\[
\mathcal{X} \rightarrow  \Spe(A)
\]
parametrized by the base scheme $\Spe(A)$, together with a group homomorphism $G\rightarrow \A _A(\mathcal{X})$ such that there is a 
$G$-equivariant isomorphism $\phi$ 
from 
the fibre over the closed 
point of $A$ to the original curve $X$:
\[
\phi: \mathcal{X}\otimes_{\Spe(A)} \Spe(k)\rightarrow X. 
\]
Two deformations $\mathcal{X}_1,\mathcal{X}_2$ are considered 
to be equivalent if there is a $G$-equivariant isomorphism $\psi$, 
making the following diagram  commutative:
\[
\xymatrix{
\mathcal{X}_1 \ar[rr]^{\psi} \ar[dr] & & \mathcal{X}_2 \ar[dl] \\
& \Spe A &
}
\]

Let $\mathcal{C}$ denote the category of local Artin algebras over $k$.
The  global deformation functor is defined:
\[
\Dgl: \mathcal{C} \rightarrow \rm{Sets}, 
A \mapsto
\left\{
\mbox{
\begin{tabular}{l}
Equivalence classes \\
of deformations of \\
couples $(X,G)$ over $A$
\end{tabular}
}
\right\}
\]

The deformation functor $\Dgl$ of non singular curves together with a 
subgroup of the automorphism group, admits a pro-representable  
hull $R$ as J. Bertin and A. M\'ezard \cite{Be-Me}  proved 
using Schlessinger's  \cite{Sch} approach.

The Krull dimension of the hull is in general smaller than the 
dimension of the tangent space of the deformation functor;
there are obstructions preventing infinitesimal deformations to be 
lifted. 

The study of obstructions, and therefore the study of the ring 
structure of the hull, is a difficult question but it has been done 
for ordinary curves by G.Cornelissen and F. Kato \cite{CK}.

J. Bertin and A. M\'ezard proved a local-global 
principle \cite[th. 3.3.4,cor. 3.3.5]{Be-Me}, that  reduces the problem 
of  computing the Krull dimension, to the computation of the dimension of 
a local deformation functor attached to each wild ramified point. 
The main difficulty for the study 
of obstructions, is that the  obstructions to such liftings are 
reduced to  obstructions  of lifting representations from  
the groups $\A k[[t]]$ to $\A A[[t]]$, where $A$ is a local Artin algebra
with  maximal ideal $m$, such that $A/m=k$.
The automorphism 
group of the ring of formal powerseries is a difficult object to study. 
In this paper we try to reduce the problem of these liftings to a similar lifting 
problem involving general linear groups.

In order to make the presentation and the calculations simpler, 
we will restrict ourselves
to subgroups of the automorphism group that 
have order a power of $p$. Equivalently, for all  decomposition groups at various 
points $P$ of the curve,
we  assume $G_0(P)=G_1(P)$, where $G_i(P)$ denotes the ramification filtration 
of the decomposition group at a ramified point.

We will give criteria depending on the Weierstrass semigroup 
structure on the wild ramified points, (Prop. \ref{semigroup-criter},
Cor.
\ref{semigroup-criter2}) that allow us to connect  
liftings of representations to $\A A[[t]]$ with 
 liftings of representations of  general linear groups and 
the deformation theory of such representations is better
understood. 
If these conditions do not hold, then we can restrict ourselves to a subfunctor 
$\mathcal{D} \subset D$ by posing the desired conditions as deformation conditions
in the sense of B. Mazur \cite[p.289]{MazDef}. Then the new subfunctor also has a hull $R_{\mathcal{D}}$ and 
$\dim_k R_{\mathcal{D}} \leq \dim_k R$.

We have to mention that 
our  approach does not give us the ring structure of the hull, 
but provides  us with a method   to compute its Krull dimension. 

The Krull dimension has been considered, by different methods, by R. Pries \cite{Pries:04},
but only for deformations that do not split the branch locus. The conditions we put 
are more general than the ones of Pries  and  can be applied to a wider class of deformations.

We begin our study by showing that there is a faithful 
representation 
\[
\rho:G_1(P) \rightarrow GL_n(k),
\]
of the $p$-part of the decomposition 
group at a wild ramified point to a suitable general linear group, 
and we show how to relate the filtration of the ramification 
group to the radical decomposition of the algebra of 
lower triangular matrices. 
This allows us to describe the structure of $G_1(P)$ for some 
interesting examples. In particular we are able to prove that 
the representation in the case of ordinary curves is 
three dimensional.

Next we give conditions so that deformations over a local  domain $A$, respect 
the flag of the Weierstrass subspaces and give rise 
to lifting of $\rho$ to a representation
\[
\tilde{\rho}:G_1(P) \rightarrow GL_n(A),
\]
that reduces to $\rho$ modulo $m$. 

In order to perform such a construction we need deformations
over local domains that can be extended to the generic fibre.
 Unobstructed deformations give rise to families 
over formal schemes that do not possess a generic fibre. We employ 
Artin's algebraisation theorem in order to show that an extension of the
family 
to the generic fibre is always possible. As a result, we can use 
(prop. \ref{bertin-gen}) an algebraic equivalence 
argument in order to compare Artin's representations at wild 
ramified points in both the generic and the special fibre, 
obtaining a generalisation of a theorem of Bertin.

We introduce a deformation functor $F(\cdot)$ for deformations 
of matrix representations and we are able to relate 
these two deformation functors in proposition
\ref{relateFD}.  

The deformation ring  $R_{\mathcal{D}}$ corresponding to the 
subfunctor $\mathcal{D}$  may have nilpotents
\cite[4.4.1]{CK}
and in general might not be irreducible.
Factoring out the radical of $R_{\mathcal{D}}$,
we obtain a finitely generated $k$-algebra
without nilpotens that corresponds to 
finite union of irreducible sets.
In proposition  \ref{prop123} and 
corollary \ref{cor123} we prove 
that all these algebraic sets have equal 
dimension, and this dimension is  equal to the dimension of the  
tangent deformation functor 
$F(k[\epsilon])$.
Moreover, we prove that  infinitesimal deformations in 
$F(k[\epsilon])$ are unobstructed and this fact  allows 
us to compute the desired Krull dimension.

Many authors have tried to study the ``simplest possible'' wild ramification.
From the point of view of ramification filtrations the simplest wild ramification 
is the ``moderate'' wild ramification, {\em i.e.}, when for all wild ramified 
points $P$ we have $G_i(P)=\{1\}$ for all $i\geq 2$. For instance 
ordinary curves have this property. 

From the representation point of view the simplest possible wild ramification at a point $P$ is 
when the faithfull representation  attached to this point is two  dimensional. 
We can prove that this is the case for all curves of genus $g$ so that $p> 2g-2$.
One of the new results of this article is the following:
\begin{theorem}
If the faithfull representation attached to a wild ramification point is two dimensional,
then the hull of the local deformation functor at this point is the ring of formal power series 
$k[[x_1,\ldots, x_s]]$, where $s=\log_p |G_1(P)|$. 
\end{theorem}

We also try to illustrate our method by giving examples and by comparing our 
method with computation done by other authors.  Namely we  apply our method to the case 
of deformations
of ordinary curves and 
as a final application, 
we show how the tools we have developed can be used in order to study deformations of the curves 
$y^p-y=\sum_{i=1}^{s-1}x^{p^i+1} +x^{p^s+1}$.

{\bf Acknowledgments} The author would like  to thank the participants of the conference
in Leiden on 
{\em Automorphisms of Curves} for enlightening 
conversations and especially R. Pries and M. Matignon for their corrections and remarks.

\section{Representations}
Let $C$ be a nonsingular complete curve defined over an algebraically 
closed field of characteristic $p\geq 0$. 
Let $G$ be a subgroup of the automorphism group of $C$, 
and let $P$ be a wildly ramified point of $C$. 
We denote by $G_0(P)$ the decomposition 
subgroup of $G$, and by  $G_1(P)$   the $p$-part of $G_0(P)$. 
The $p$-part of the decomposition 
group can be analysed in terms of the sequence 
of the $i$-th ramification groups \cite[chap. IV]{SeL}:
\begin{equation} \label{ram-fil-seq}
G_1(P) \geq G_2(P) \geq \cdots \geq G_n(P)  > \{1 \}.
\end{equation}

For every point $P$  of the curve $C$ of genus $g$ 
we consider the sequence of $k$-vector 
spaces 
\begin{equation} \label{flag}
k=L(0)=L(P)= \cdots  =L((i-1)P) < L( iP) \leq  \cdots \leq  L( (2g-1)P ),
\end{equation}
where
\[
L(iP):=\{ f \in k(C): {\rm div}(f) +iP \geq 0 \}=H^0(X,\cL(iP)).
\]
It is known that there are exactly $g$ pole numbers that are smaller or equal to $2g-1$.  If $g\geq 2$ then there is at least 
one of them not divisible by the characteristic.
\begin{lemma} \label{fai-rep-lemma}
        Let $m\leq 2g-1$ be the smallest  pole number not divisible by the characteristic. 
There is  a faithful representation 
\begin{equation} \label{fai-rep}
\rho: G_1(P) \rightarrow \mathrm{GL} \big( L(mP) \big)
\end{equation}
\end{lemma}
\begin{proof}
It is clear that the space ${L}(mP)$ is preserved 
by any automorphism in $G(P)=G_1(P)$.  Let $f$ be a function with pole at $P$ of order $m$. 
We can write $f$ as 
$f=u\frac{1}{t^m}$, where $u$ is a unit in the local ring $\O_P$. 
Since $(m,p)=1$, Hensel's lemma implies that $u$ is an $m$-th power
so the  
local uniformizer  can be selected  so that  $f=\frac{1}{t^m}$. 
Let $\sigma \in G_1(P)$ such that $\sigma(1/t^m)=1/t^m$. 
Then $\sigma(t)= \zeta t$ where $\zeta$ is an $m$-th root of unity. 
Therefore, if $\sigma$ induces the trivial matrix 
in $\A {L}(mP)$  and $\sigma$ is of order $p$, then $\zeta=1$ 
since $(p,m)=1$. 
\end{proof}
The above lemma makes the $p$-part of the decomposition 
group $G(P)$ realizable as a finite  algebraic subgroup of the linear  group 
$GL_n(k)$. 
Moreover the flag of vector spaces $L(iP)$ for $i\leq m$ is preserved, 
so the representation matrices are upper triangular. 

We assume that 
 $m=m_0>m_1>\cdots >m_r=0$, 
are the pole numbers less than $m$. Therefore, a basis for 
the vector space $L(mP)$ is given by 
\[\{1,u_i\frac{1}{t^{m_i}},\frac{1}{t^m}:\mbox{ where } 1<i<r, p \mid m_i \mbox{ and } u_i 
\mbox{ are  units} \}\]
According to this basis, an element  $\sigma \in G_1(P)$ acts on $L(mP)$ by 
\[
\sigma\frac{1}{t^m}=\frac{1}{t^m}+ \sum _{i=1}^r c_i(\sigma) 
u_i \frac{1}{t^{m_i}}, 
\]
and equivalently it maps the local uniformizer $t$ to 
\[
\sigma(t)= \frac{\zeta t}{\left( 1+\sum_{i=1}^r c_i(\sigma)
u_i t^{m-m_i} \right)^{1/m},  } 
\]
where $\zeta$ is an $m$-th root of one.

The  above expression can be written in terms 
of a formal power series as:
\begin{equation} \label{eqsform}
\sigma(t)=\zeta t\big( 1+ \sum_{\nu \geq 1} a_\nu(\sigma) t^{\nu} \big).
\end{equation}
On the other hand the composition of the formal powerseries 
\[
f= t \sum_{i\geq 0} a_i t^i \mbox{ and } g = t \sum_{j \geq 0} b_j t^j 
\]
is written as $t a_0 b_0 + \cdots $, so the automorphism $\sigma ^{|G_1(P)|}=1$ 
as it is 
given in (\ref{eqsform}) is $t \zeta^{p^l}+ \cdots = t$ and since $m$ is 
prime to $p$, $\zeta=1$ and (\ref{eqsform}) can be written as
\begin{equation} \label{eqsform2}
\sigma(t)= t\big( 1+ \sum_{\nu \geq 1} a_\nu(\sigma) t^{\nu} \big).
\end{equation}

The above computation allows us to compute the ``gaps" in the filtration 
of the group $G_1(P)$. 
\begin{pro} \label{gap-comp}
Let $P$ be a point on the curve $C$ and let 
 \[\rho: G_1(P) \rightarrow GL_{\dim(L(mP))} (k)\] be the corresponding  faithful 
representation we considered  on lemma \ref{fai-rep-lemma}. 
Let $m=m_0>m_1 > \cdots > m_r=0$ be the pole numbers at $P$ that are 
smaller than $m$. For the filtration of $G_1(P)$ we have
\[
\mbox{if }G_{i}(P)  >  G_{i+1}(P) \mbox{ then  } i+1=m-m_k+1, 
\]
for some pole number $m_k$. 
\end{pro}
\begin{proof}
By definition $\sigma \in G_i(P)$ \cite[prop. 1 p.62]{SeL} if and only if
$\sigma (t)-t\in t^{i+1} k[[t]]$.
Notice that there is at least one $c_i(\sigma)\neq 0$, because if all 
$c_i(\sigma)=0$ for $i=1,...,r$ then $\sigma(1/t^m)=1/t^m$  
and $\sigma$ is the identity. 
The valuation of the  expression $\sigma(t)-t$ can be explicitly 
computed:
\[
\sigma(t)-t=-\frac{1}{m} \sum_{i=1}^r c_i(\sigma) u_i 
t^{m-m_i+1}+\cdots ,\]
therefore
\[
v_t \big (\sigma(t)-t \big)=m-m_k+1
\]
where $k=\min\{ i : c_i(\sigma) \neq 0\}$.
 The possible valuations are given by:          
\[
\begin{array}{ll}
m-m_1+1 & \mbox{if } c_1(\sigma)  \neq 0 \\
m-m_2+1 & \mbox{if } c_2(\sigma) \neq 0, c_i(\sigma)=0 \mbox{ for } i<2 \\
\vdots & \vdots \\
m -m_r+1 & \mbox{if } c_r(\sigma) \neq 0,  c_i(\sigma)=0 \mbox{ for } i<r
\end{array}
\]
Assume that $\sigma \in G_{i}(P)$ but $\sigma \not \in G_{i+1}(P)$, 
thus $v(\sigma(t)-t)=i+1$ and this equals some $m-m_k+1$. 
\end{proof}
\begin{cor} \label{nodivart}
Every jump $i$  in the  
 ramification filtration, {\em i.e.},
$G_i > G_{i+1}$ is not divisible by $p$. 
\end{cor}
\begin{proof}
By lemma \ref{gap-comp} every gap in the ramification 
filtration is given as $m-m_k$, where $m$ is not divisible by $p$
and $m_k$ are divisible by $p$ \cite[IV. prop. 11]{SeL}.
\end{proof}

{\bf Examples:} 
\\
{\bf 1.} \label{FermatDEF} The Fermat curves $x^n+y^n+1=0$, 
where $n-1=p^h$. The automorphisms of these curves where studied by  H. W. Leopoldt in 
\cite{Leopoldt:96}, even if $n-1$ is not a power of the characteristic. Leopoldt  
constructed a basis for the space of holomorphic differentials of the curve and 
he was able to prove that for the points of the form $P:(x,y)=(\zeta_{2n},0)$
where $\zeta_{2n}$ is a $2n$-root of one, we have the following sequence 
of $k$-vector spaces \cite[Satz 4]{Leopoldt:96}:
\[
k={L}(0P)={L}(P)=\cdots={L}( (n-2)P)<
{L}( (n-1)P) < {L} (n P) \leq \cdots
\]
The interesting case for us (Hermitian Function Fields) appears when 
$n-1$ is a power of the characteristic, so in this case the representation 
of the decomposition subgroup is of the form:
\[
\rho: G_0(P) \rightarrow GL\big( \mathcal{L} (nP ) \big)
\]
with 
\[
\sigma \mapsto \left( 
\begin{array}{ccc}
1   & 0 & 0 \\
\alpha & \chi & 0 \\
\gamma & \beta & \psi
\end{array}
\right).
\]
According to proposition \ref{gap-comp} the filtration of the decomposition 
subgroup is given by:
\[
G_0(P) > G_1(P) > G_2(P) = \cdots =G_n(P) >  G_{n+1}(P) = \{1\}, 
\]
{\em i.e.}, the gaps of the filtration are in $0,1,n$. 
\\
{\bf 2.}\label{n2} The curves $x^n+x^m+1=0$, where $m \mid n$ and $m-1=p^h$
The automorphism group of a nonsingular model of the above curve
is studied by the author, in 
\cite{Ko:98}. It is proved that for the points $P: (x,y)=(\zeta_{2n},0)$
we have the following sequence of vector spaces \cite[eq. (4)]{Ko:98}
\[
k={L}(0P)={L}(P)=\cdots = {L} ( (m-1) P ) < 
{L} (m P) = {L} ((m+1)P) \leq \cdots
\]
Since $m$ is not divisible by $p$ we have the following representation
\[
G_0(P) \rightarrow GL\big( {L} ( mP ) \big), 
\]
sending 
\[
\sigma \mapsto  \left(
\begin{array}{cc}
1 & 0 \\
\alpha & \chi
\end{array}
\right).
\]
Thus $G_0(P)$ is the semidirect product of an elementary abelian group 
by a cyclic group of order prime to the characteristic. 
For the ramification filtration of $G_0(P)$ we have
\[
G_0(P) > G_1(P) > G_2(P) = \cdots =G_{m} (P)  > \{1 \}. 
\]
\\
{\bf 3.} Ordinary Curves. \label{OrdCurves} 
A curve is called ordinary if the $p$-rank of the Jacobian is 
equal to the genus of the curve. 
It is known that ordinary curves form a Zariski-dense set 
in the moduli space of curves of genus $g$. 
For ordinary curves we have that $G_2(P)=\{1\}$ \cite{Nak}, thus
we have the following picture for the faithful representation  $\rho$
of the group $G_1(P)$: 
There is a gap at $G_1(P)>G_2(P)=\{1\}$, thus $m_i=m-1$ for some 
$i$, and this $i=1$. In other words the pole numbers that are smaller or 
equal to $m$ are $\{m,m-1\}$. 
This implies that if the genus $g$ of $X$ is $g\geq 1$ 
then
the representation has dimension  3, 
because otherwise,{\em i.e.}, if the representation is two 
dimensional, we have the following sequence 
\[
k=L(0P)=L(P)=\cdots L( (m-1)P) < L( (m)P). 
\]
But $m-1$ is a pole number so $m-1=0$ and 
$m=1$, {\em i.e.}, the Weierstrass semigroup 
is the semigroup of natural numbers, a contradiction, for  $g\geq 2$. 

Moreover if $c_1(\sigma)=0$ then all $c_i(\sigma)=0$ for $i>1$, otherwise 
there  will be  more jumps at higher groups $G_i$, and this is impossible.
This proves  that $c_1(\sigma)=0$ if and only if $\sigma=0$. By multiplying 
the representation matrices we can easily deduce that the map 
\begin{equation} \label{c1}
c_1: G_1(P) \rightarrow k,
\end{equation}
is a faithful homomorphism  of the elementary abelian group $G_1(P)$ into the 
additive group of $k$. 

The representation matrices are commuting, and by computation 
this implies that all elements  $\rho_{j+1,j}(\sigma)$, 
of the representation matrix, are of the 
form $\lambda_j c_1(\sigma)$, and  $\lambda_j$ is independent of 
$\sigma$. 

Since  $c_1$ is faithful character, such that
$c_1(\sigma)=0$ implies that 
$\rho_{ij}(\sigma)=0$, for all $i\neq j$, 
we can write 
$\rho_{ij}(\sigma)=c_1(\sigma) a_{ij}(\sigma)$.

{\bf 4.} \label{p-cycl} p-cyclic covers of the affine line. 
In this example we apply our computations to Artin-Schreier curves that 
are nonsingular models of the function field defined by:
\[C_{t_1,\ldots t_{s-1}}:W^p-W= \sum_{i=1}^{s-1} t_i X^{p^i+1} +X^{p^s+1}.\] 
These curves give extreme examples of automorphism groups and 
were studied by H. Stichtenoth  \cite{StiII} and C. Lehr, M.Matignon \cite{CL-MM}, N. Elkies \cite{Elkies99}, 
van der Geer and van der Vlught \cite{Geer-Vlught}.

There is only one ramified point in the cover $C_{t_1,\ldots t_{s-1}} \rightarrow 
\mathbb{P}^1$, the point $P$ that is over the point $X=\infty$ of $\mathbb{P}^1$. 
The Weierstass semigroup at $P$ is computed by H. Stichtenoth \cite{StiII} to 
be equal to $(p^s+1) \mathbb{N} + p \mathbb{N}$.

Thus, the smaller pole number that is not divisible by $p$ is $p^s+1$ and 
the Weierstrass semigroup up to $p^s+1$ is computed to be
\[
0,p,2p,\ldots, \left[ \frac{1+p^s}{p}\right]p, 1+p^s.
\]
One can prove that $\left[ \frac{1+p^s}{p}\right]p=p^s$.
The remainder of the division of $p^s+1$ by $p$, 
is $1$.
According to proposition \ref{gap-comp} the possible gaps at the ramification 
filtration are at the numbers $1+kp$, 
$k=0,\ldots,\left[ \frac{1+p^s}{p}\right]=p^{s-1}$. 
The dimension $\dim_k L \big((p^s+1)P \big)$ is  $n=\left[ \frac{1+p^s}{p}\right]+2=p^{s-1}+2$
and 
the representation of $G_1(P)$ to $L\big( (p^s+1) P\big)$ is given  
by an $n\times n$ lower triangular matrix with $1$ in the diagonal.

More precisely,  if we choose the natural basis $\{1,X,X^2,\cdots,X^{p^{s-1}},W\}$
of 
$L \big((p^s+1)P \big)$ then the representation $\rho$ is given by the matrix
\begin{equation} \label{repCM}
	\rho(\sigma)_{ij}=
	\left\{
	\begin{array}{ll}
		0 &  \mbox{ if } i<j \\
		1 &  \mbox{ if } i=j \\
		y^j \binom{i}{j} & \mbox{ if }  i>j, i\neq p^{s-1}+1 \\
		b_j(y)  & \mbox{ if }  i=p^{s-1}+1>j
	\end{array}
	\right.,
\end{equation}
where $b_j(y)$ are the coefficients of the polynomial $P_{f }(X,y)$,
and $y$ is a solution of $Ad_{f}(Y)=0$ as defined in lemma 4.1 and definition 
4.2 in \cite{CL-MM}.

%
%
%
\section{Deforming Branch Points}
Assume that we have a deformation $\mathcal{X} \rightarrow \Spe A$ over an 
Artin local ring, that admits a fibrewise  action of the  $G\subset \A (X)$.
In lemma \ref{fai-rep-lemma}. we have assigned  to every wild  ramified point 
$P$  a representation  of  the decomposition group $G(P)$ that corresponds to an upper triangular matrix.
We will try to lift this representation to a representation of an upper triangular matrix with coefficients in $A$.
 This will be accomplished by proving that we can 
fullfill the requirements of the following 
\begin{pro}\label{crit-matrix-lift}
        Consider a deformation $\mathcal{X} \rightarrow \Spe A$ of the curve $X$ 
        defined over a local integral domain $A$, with $G \subset \A(X)$ acting 
        fibrewise on $\mathcal{X}$, and let $P$ be a wild ramified point on the 
        special fibre of $\mathcal{X}$. 
	Suppose that there is a sequence of $G$-invariant
invertible  $\O_{\mathcal{X}}$-modules  $\mathcal{L}_i$ so that the corresponding spaces of global sections $L_i:=\Gamma(\mathcal{X},\mathcal{L}_i)$
satisfy:
        \[
        L_1 \subset L_2 \subset \cdots \subset L_n,
        \]
        and  $L_i$ are finitely generated free $A$-modules and $L_i \otimes_A k= L(iP)$. 
Then the faithful  representation defined in (\ref{fai-rep}), 
can be lifted to a representation:
\[
\rho_1:G_1(P) \rightarrow GL_n(A),
\]
such that $\rho_1(\sigma), \sigma \in G_1(P)$ is a lower triangular matrix with $1$ at 
the diagonal.
\end{pro}
\begin{proof}
	Let $k=A/m$, be the residue field of the local ring $A$, modulo 
	the maximal ideal $m$ of $A$. 
	The $A$-module $L_i$ is  a finitely generated free $A$-module,
	therefore 
	$\dim_k (L_i \otimes_A k)=\mathrm{rank}_A L_i$.
        Moreover, since $L_i$ are $G$-modules, there is a natural representation 
	$\rho:G_1(P) \rightarrow GL_n(A)$. 
Since the flag $L_i$ of $A$ modules is preserved, the 
representation matrices of $\rho_1$ are lower 
triangular, and elements on  the diagonal must 
have a multiplicative $p$-group structure, so they are units.
\end{proof}

In subsection \ref{sub-sec1} we recall the theory of effective Cartier divisors
and we define the relative ramification divisor. We also show that if exists a deformation  $\mathcal{X} \rightarrow \Spe A$ over an integral domain $A$ 
  that admits a fibrewise $G$-action exists
then there is a constrain in the Artin representations at the special and the generic fibres expressed in proposition 
\ref{bertin-gen}.

Subsection \ref{sub-sec2} is devoted to the problem of algebraisation. In order to construct the desired representation 
we will need deformations that possess generic fibres. The passage from deformations defined over formal schemes 
to deformations that defined over ordinary schemes is given by Artin's algebraisation theorem \ref{ACOA} and 
in subsection \ref{sub-sec2} we check that Artin's algebraisation theory can be applied.

In subsection \ref{sub-sec3} we try to construct the sequence $L_i$ of invertible  $\O_{\mathcal{X}}$ modules of proposition 
\ref{crit-matrix-lift}. For this we need Grauert theorem \cite[III.12.9]{Hartshorne:77} and this is the reason we want our deformations to possess 
a generic fibre.    We give a criterion \ref{semigroup-criter}  for the construction of the sequence $L_i$ and  some lemmata 
that imply the truth of this criterion  and depend only of the form of the ramification filtration of the special fibre.
 Unfortunately it seems that there are bad deformations and for them 
the construction of proposition \ref{crit-matrix-lift} is not possible.  The best think we can do is to show that 
the criterion \ref{semigroup-criter} defines a deformation condition in the sense of Mazur, which in turn gives rise to a 
subfunctor that can be handled with our tools. Of course, this approach will lead to a  lower bound 
of the desired Krull dimension. 

\subsection{Relative Ramification divisors} \label{sub-sec1}
In this section we would like to address the following question:
Let $P$ be a wild ramified point, of the special fibre of the
deformation in study, with decomposition group
$G(P)$. Is it possible to find a ``horizontal 
 divisor"  $\bar{P}$ that is invariant under the action of
$G(P)$, so that the intersection of $\bar{P}$ with the special fibre is the original point $P$?

Let $\mathcal{X} \rightarrow \Spe A$ be an $A$-curve, admitting a 
fibrewise action of the finite group $G$, where $A$ is a 
Noetherian local ring. 
Let $\mathcal{Y} \rightarrow \Spe A$ be the quotient $\Spe A$-curve.   
A good notion of a horizontal divisor in this case can be given 
in terms of effective Cartier divisors.
An effective Cartier divisor $D$ on $\mathcal{X}\rightarrow \Spe A$,
is a closed subscheme $D\subset \mathcal{X}$, such that $D$ is flat 
over $\Spe A$, and the ideal sheaf $I(D)\subset \O_{\mathcal{X} }$ is 
an invertible $\O_{\mathcal{X}}$-module. 

We would like to assign to the cover of $A$-schemes $\mathcal{X} \rightarrow 
\mathcal{Y}$ a ramification divisor $D_{\mathcal{X}/\mathcal{Y}}$, such that the intersection 
of 
$D_{\mathcal{X}/\mathcal{Y}}$
 with the fibres of the morphism 
$f: \mathcal{X} \rightarrow \mathcal{Y}$, 
corresponds to the usual notion of ramification divisor for coverings 
of $k$-curves. 

Let $S=\Spe A$, and $\Omega_{\mathcal{X}/S}$, $\Omega_{\mathcal{Y}/S}$ be
the sheaves of relative differentials of $\mathcal{X}$ over $S$ and 
$\mathcal{Y}$ over $S$, respectively. 

We begin by defining the ideal sheaf of the ramification divisor 
$D_{\mathcal{X}/\mathcal{Y}}$ as  
\[
\mathcal{L}(-R)= \Omega_{\mathcal{X}/S} ^{-1}
\otimes_S f^* \Omega_{\mathcal{Y}/S}. 
\]
We will first prove that 
 the above defined divisor $D_{\mathcal{X}/\mathcal{Y}}$ is indeed 
an effective Cartier divisor. 

We will use the following 
\begin{crit} \label{criter-eff-cart}
Let $S$ be locally Noetherian, $\mathcal{X}$ a flat $S$-scheme of finite type, 
and $D\subseteq \mathcal{X}$, a closed subscheme which is flat over $S$. 
Then $D$ is an effective Cartier divisor in $\mathcal{X}/S$  if and only 
if, for all geometric points $\Spe k \rightarrow S$ of $S$, the 
closed subscheme $D\otimes_S k$ of $\mathcal{X} \otimes_S k$ is an effective 
Cartier divisor in $\mathcal{X} \otimes _S k/k$. 
\end{crit}
\begin{proof}
\cite[Cor. 1.1.5.2]{KaMa}
\end{proof}


\begin{pro}
The ramification divisor $D_{\mathcal{X}/\mathcal{Y}}$ is an 
effective $A$-divisor. 
\end{pro}
\begin{proof}

We are interested in deformations of nonsingular curves. 
Since the base is a local ring and the special fibre is nonsingular, 
the deformation $\mathcal{X} \rightarrow \Spe A$ is smooth.  
(See the remark after the definition 3.35 p.142 in \cite{LiuBook}).
The smoothness of the curves $\mathcal{X}\rightarrow S$, 
and $\mathcal{Y}\rightarrow S$, implies that the sheaves 
$\Omega_{\mathcal{X}/S}$ and $\Omega_{\mathcal{X}/S}$ are $S$-flat,
\cite[cor. 2.6 p.222]{LiuBook}. 

On the other hand the sheaf $\Omega_{\mathcal{Y},\Spe A}$ is 
by \cite[Prop. 1.1.5.1]{KaMa}  $\O_{\mathcal{Y}}$-flat. 
Thus, $\pi^*(\Omega_{\mathcal{Y}, \Spe A})$ is $\O_{\mathcal{X}}$-flat
and therefore $\Spe A$-flat \cite[Prop. 9.2]{Hartshorne:77}.
The desired result follows by 
 criterion \ref{criter-eff-cart}.
\end{proof}


%
%
%
%
%
%
Let $A$ be a local $k$ algebra, that is a domain,  and let 
$\mathcal{X} \rightarrow \Spe A$ be a deformation of the curve
$X$, and  $P\in X$, be a wild ramified point with decomposition group 
$G_0(P)$. Assume that the point $P$ appears in the ramification divisor 
of the covering of the special fibres 
$X\rightarrow X/G$, with multiplicity  $d$, given by Hilbert's formula 
\cite[III.8.8]{StiBo}.

We consider the effective Cartier divisor
$D_P=\sum_{j=1}^\lambda a_i \bar{P}_i$,  
where $\bar{P}_i$ denotes the irreducible 
components of the ramification 
divisor $R=R_{\mathcal{X}/\mathcal{Y}}$ that
intersect the special fibre of $\mathcal{X}$ at 
$P$.

We have the following picture:

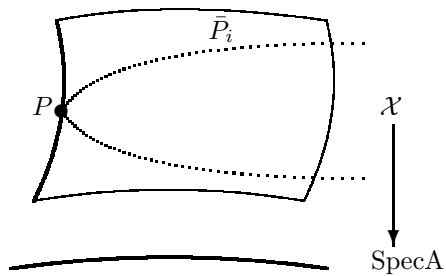
\setlength{\unitlength}{0.6mm}
\begin{picture}(70,90)(-50,-30)
\linethickness{1pt}
\qbezier(0,0)(10,20)(5,40)
\qbezier(-5,-15)(30,-10)(65,-15)
\thinlines
\qbezier(60,0)(70,20)(65,40)
\qbezier(5,40)(30,45)(65,40)
\qbezier(0,0)(30,5)(60,0)
\thicklines
\put(6.37,20){\circle*{3} }
\put(-2,19){ $P$ }
\put(75,19){ $ \mathcal{X}$ }
\put(80,17){\vector(0,-1){27} }
\put(75,-15) {${\rm Spec A}$} 
\qbezier[50](6.37,20)(15,35)(73,35)
\qbezier[50](6.37,20)(15,5)(73,5)
\put(39,36){$\bar{P}_i$ }
\end{picture}

Two horizontal branch divisors can collapse to the same point in 
the special fibre. For instance this always happens if a 
deformation of  curves from positive characteristic to characteristic zero
with a wild ramification point is possible.

For a curve $X$ and a branch point $P$ of $X$ we will 
denote by  $i_{G,P}$  the order  function of the filtration of $G$ at $P$. 
The integer  $i_{G,P}(\sigma)$ 
is equal to the multiplicity of $P\times P$ in the intersection of 
$\Delta .\Gamma_\sigma$ in the relative $A$-surface $X \times X$, 
where $\Delta$ is the 
diagonal and $\Gamma_\sigma$ is the graph of $\sigma$ \cite[p. 105]{SeL}. 
Using an algebraic equivalence argument on 
$\mathcal{X}\times_{\Spe A} \mathcal{X}$ we see the following generalisation 
of the result of J. Bertin \cite{BertinCRAS}:
\begin{pro}\label{bertin-gen}
	Assume that $A$ is an integral domain, and let $\mathcal{X}\rightarrow \Spe A$
	be a deformation of $X$. 
        Let $\bar{P}_i$, $i=1,\cdots,s$ be the  horizontal branch divisors 
        that intersect at the special fibre, at point $P$, and let $P_{i}$ be 
        the corresponding points on the generic fibre. For the Artin 
        representations attached to the points $P,P_{i}$ we have:
        \[
        \mathrm{ar}_P(\sigma)=\sum_{i=1}^s \mathrm{ar}_{P_{i}}(\sigma). 
        \]
\end{pro}

{\bf Remark} \label{ex-cor} Consider the case of deformations of
ordinary curves, together 
with a $p$-subgroup of the group of automorphisms. 
Then $|\mathrm{ar}_P(\sigma)|=2$  for all 
$\sigma \in G(P)=G_1(P), \sigma\neq 1$ \cite{Nak}.
On the other hand the ramification at the 
points of the generic fibre is also wild and \ref{bertin-gen} implies that 
there is only one horizontal branch divisor extending every wild ramification 
point $P$. 
%
%
%

%
%

\subsection{Application of Algebraisation Theory} \label{sub-sec2}
%
%

We would like to have our curves deformed 
over bases that have generic fibres. Formal deformation theory 
gives us information whether a curve can be defined over the formal 
spectrum of a complete ring, {\em i.e.} over 
\[
\mathrm{Spf} R=\{ P \in \Spe R: P \mbox{ is open with respect to the } m_R-adic \mbox{ topology} \}, 
\]
where $R$ is a complete domain with maximal ideal $m_R$. In order to extend the 
family over the generic fibre (the $0$ ideal is not open) we have to use 
an algebraisation argument. 
 
Let us denote by $\mathcal{C}$ the category of local Artin $k$-algebras
and by $\hat{\mathcal{C}}$ the category of complete Noetherian local 
$k$-algebras. 
A covariant functor $F:\mathcal{C}\rightarrow \mathrm{Sets}$ can 
be extended to a functor $\hat{F}:\hat{\mathcal{C}}\rightarrow \mathrm{Sets}$
by defining 
\[
\hat{F}(R)= \lim_{\leftarrow} F\left(\frac{R}{m_R^{n+1}}\right),
\]
where $m_R$ is the maximal ideal of $R$ and $\frac{R}{m_R^{n+1}}$
is a local Artin $k$-algebra.
On the other hand a functor $F:\hat{\mathcal{C}}\rightarrow \mathrm{Sets}$ 
induces by reduction a functor $F\mid_{\mathcal{C}}:\mathcal{C} \rightarrow \mathrm{Sets}$.
For any covariant functor $F:\hat{\mathcal{C}}\rightarrow \mathrm{Sets}$
there is a canonical map 
\[
F(R) \rightarrow \hat{F}(R)=\lim_{\leftarrow} F\left(\frac{R}{m_R^{n+1}}\right).
\]
The above map is not in general a bijection. Let us denote by  
$h_R(\cdot)=Hom(R,\cdot)$.
One can also prove \cite[lemma 2.3]{Kai-wen-lan} that  
$\hat{F}(R)\stackrel{\cong}{\rightarrow} Hom(h_R,F)$.
A functor $F:\mathcal{C}\rightarrow \mathrm{Sets}$ is called prorepresentable 
if there is an $R\in Ob(\mathcal{C})$ and  $\hat{\xi} \in \hat{F}(R)$, that
induces an isomorphism 
\[
\hat{\xi}: h_R(A)\stackrel{\cong}{\rightarrow} F(A).
\]
A formal deformation of $\xi_0 \in F(k)$ is an element 
$\hat{\xi}\in \hat{F}(R)$, where $R\in Ob(\hat{\mathcal{C}})$.

Let $F$ be a functor $F:\hat{\mathcal{C}}\rightarrow \mathrm{Sets}$. We 
say that $F\mid_{\mathcal{C}}$ is effectively prorepresentable if 
$F\mid_{\mathcal{C}}$ is prorepresentable by $\xi \in \hat{F}(R)$, and this 
$\xi$ is the image of an element in $F(R)$, under the map 
$F(R)\rightarrow \hat{F}(R)$.

J. Bertin, A. M\'ezard \cite{Be-Me}, 
introduced at the wild ramified point $P$ the deformation functor:
\begin{equation} \label{Bertin-Mezard-functor}
D:\mathcal{C} \rightarrow {\rm Sets}, 
A \mapsto 
\left\{
\mbox{
\begin{tabular}{l}
lifts $G\rightarrow \A (A[[t]])$ of $\rho$ mod- \\ulo 
conjugation with an element \\ of $\ker(\A A[[t]]\rightarrow k[[t]] )$
\end{tabular}
}
\right\}
\end{equation} 
\begin{lemma} \label{can-map1}
Let $D$ be the functor  defined in (\ref{Bertin-Mezard-functor}). For every 
complete $k$-algebra $A$ the canonical map 
\begin{equation} \label{can-map}
D(A) \rightarrow \hat{D}(A).
\end{equation}
is a bijection. 
\end{lemma}
\begin{proof}
Let $A$ be a complete $k$-local algebra with maximal ideal $m_A$ and 
denote by $A_n$ the quotient $A/m_A^n$. 
We will show first that the map in (\ref{can-map}) is surjective.  
Let   $\rho_n: G \rightarrow \A A_n[[t]]$ be a system of representatives of maps  
such that the following diagram 
\[
\xymatrix{
 A_{n+1} \ar[r]^{\rho_{n+1}(g)}  \ar[d]_{\mod m_A^n} &  A_{n+1} \ar[d]^{\mod m_A^n} \\
 A_n \ar[r]_{\rho_n(g)}  & A_n 
}
\]
is commutative. 

In order to define $\rho: G \rightarrow \A A[[t]]$ we have to define it on $t$. 
Let us write 
\[
\rho_n(g)(t) = \sum_{i=0}^\infty a_{i,n}(g) t^i. 
\]
The elements $\{a_{i,n}(g)\}_n$ form an inverse system and give rise to 
a limit element $a_i(g)\in A= \lim_{\leftarrow} A_n$ and to the desired extension
\[
\rho(g)(t)= \sum_{i=0}^\infty a_i(g) t^i.
\]
This proves that the canonical map in (\ref{can-map}) is indeed surjective. 

In order to prove that it is also injective  we consider two representations
\[
\rho_{1},\rho_2:G \rightarrow \A A[[t]],
\]
such that for every $n$ there are isomorphisms  $\gamma_n: A_n[[t]] \rightarrow A_n[[t]]$
that induce the identity on $k[[t]]$ such that 
\[
\rho_{1,n}= \gamma_n \rho_{2,n} \gamma_n^{-1}.
\] 
Arguing in the same way as in the proof of the subjectivity we see that 
the maps $\{\gamma_n\}_n$ give rise to a well defined isomorphism  
\[
\gamma: A[[t] \rightarrow A[[t]]
\]
that induces the identity map on $k[[t]]$ and 
makes $\rho_1,\rho_2$ equivalent.
\end{proof}

\begin{orism} 
We will say that the functor $D:\hat{\mathcal{C}}\rightarrow \mathrm{Sets}$ is locally of finite
presentation if and only if for every direct limit $\lim_{\rightarrow}A_i$ of objects in $\bar{\mathcal{C}}$
the natural map 
\[
\lim_{\rightarrow} D( A_i) \rightarrow D(\lim_{\rightarrow} A_i),
\]
is an isomorphism.
\end{orism}
\begin{orism}
We will say that the functor  $D:\hat{\mathcal{C}}\rightarrow \mathrm{Sets}$ is coherent 
\cite{Hart:Coherent}
if there are two representable functors $h_2,h_1$ such that 
\[
h_2 \rightarrow h_1 \rightarrow D \rightarrow 0.
\]
\end{orism}
\begin{lemma} \label{Mpelijohn}
Every coherent functor  $D:\hat{\mathcal{C}}\rightarrow \mathrm{Sets}$ is
locally of finite presentation.
\end{lemma}
\begin{proof}
Since $D$ is coherent there is 
an  exact sequence 
\begin{equation} \label{mpeli1}
h_2\rightarrow h_1 \rightarrow D \rightarrow 0, 
\end{equation}
where $h_1,h_2$ are representable functors. Sequence (\ref{mpeli1})
implies the following commutative diagram
\[
\xymatrix{
0 \ar[r]  & Hom(D(-),\lim \limits_{\rightarrow} A_i)  \ar[r]& Hom(h_1(-),\lim \limits_{\rightarrow} A_i) \ar[r] &
Hom(h_2(-),\lim \limits_{\rightarrow} A_i) \\
0 \ar[r] & \lim \limits_{\rightarrow}  Hom(D(-), A_i) \ar[r] \ar[u]^{f} &  
 \lim \limits_{\rightarrow} Hom(h_1(-), A_i) \ar[r] \ar[u]_{f_1}  &
\lim \limits_{\rightarrow} Hom(h_2(-), A_i) \ar[u]_{f_2}
},
\]
where the last row is exact since we are considering direct limits in the category of sets. Now 
representable functors are locally of finite presentation therefore $f_1,f_2$ are isomorphisms
and the commutativity of the diagram forces $f$ to be also an isomorphism, {\em i.e.}, 
the desired result.
\end{proof}
\begin{pro} \label{coherent1}
The global deformation functor $\Dgl$ of J. Bertin A. Mezard is coherent.
\end{pro}
\begin{proof}
We will use the notation of  \cite[sect. 5]{Be-Me}. Let $\mathcal{M}_{g,n,G}$ be the functor 
classifying equivalence of triples $[\mathcal{X}/A,\phi,\theta]$, represented by the scheme $M_{g,n}^G$.
There is the following exact sequence of functors:
\[
N_{\mathrm{GL}( \Z/n\Z)}G(-) \stackrel{\alpha_1}{\rightarrow}  \mathcal{M}_{g,n,G}(-) 
\stackrel{\alpha_2}{\rightarrow} \Dgl(-) \rightarrow 0,
\]
where $N_{\mathrm{GL}( \Z/n\Z)}G$ denotes the  normaliser of $G$ in  $\mathrm{GL}( \Z/n\Z)$ 
and $N_{\mathrm{GL}( \Z/n\Z)}G(-)$ the constant group scheme with fibre the group $N_{\mathrm{GL}( \Z/n\Z)}G$.

Indeed, for a local ring $A$, and a deformation $d:\mathcal{X}\rightarrow \Spe A$  in $\Dgl(A)$,
 we can select a level structure $\phi:\mathcal{X}[n] \rightarrow  (\Z/n\Z)^{2g}_{\Spe A}$ that 
maps on  $d$ by forgetting the extra level structure. The group $G$ can be considered as a subgroup 
of $\mathrm{GL}( \Z/n\Z)$ and we can select a map $\theta:G \rightarrow \A \mathcal{X}$ such that 
\begin{equation} \label{GL-invariance}
\phi \circ \theta(\sigma)^{-1}[n]=\sigma\circ \phi.
\end{equation}
 This proves the surjectivity of $\alpha_2$.

In order to prove the exactness at the second factor 
we consider two triples $[\mathcal{X}_i/A,\phi_i,\theta_i]$, $i=1,2$
that map on $d\in \Dgl(A)$. Thus $\mathcal{X}_1\equiv \mathcal{X}_2$ and there is an element 
$a\in \mathrm{GL}( \Z/n\Z)$ such that $\phi_1=a \circ \phi_2$. But since both $\phi_1,\phi_2$ 
have to be $\mathrm{GL}( \Z/n\Z)$ invariant, {\em i.e.}, equation (\ref{GL-invariance})
has to hold we obtain that:
\[
a \circ \sigma \circ \phi_2=a \circ \phi_2 \circ \theta(\sigma)^{-1}[n]= 
\phi_1 \circ \theta(\sigma)^{-1}[n] = \sigma\circ \phi_1= \sigma \circ a \circ \phi_2  
\]
for all $\sigma \in G$. Therefore $a\in N_{\mathrm{GL}( \Z/n\Z)}G$ and the 
desired result follows.

\end{proof}
\begin{theorem}[Artin's Criterion of Algebraisation] \label{ACOA}
	Let $F:(Sch/k)^0 \rightarrow \mathrm{Sets}$ be a functor from the 
	category of the formal schemes over the algebraically closed field $k$ 
that is locally of finite presentation over $k$.
	Let $\xi_0\in F(k)$, and suppose that an effective formal deformation 
	$(R,\xi)$ exists, where $R$ is a complete Noetherian local $k$-algebra and 
	$\xi \in F(R)$. Then there is a scheme $X$ of finite type over 
	$k$ and $x\in X$ a closed point with residue field $k$, 
	and $\tilde{\xi}\in F(X)$, such that $(X,x,\tilde{\xi})$ is 
	a versal deformation of $\xi_0$, such that $\hat{\O}_{X,x})=R$, and
	\[
	\begin{array}{ccccl}
	\Spe (R) & \stackrel{\cong}{\rightarrow} & \Spe (\hat{\O}_{X,x}) & \rightarrow X \\
	\xi & \mapsto  &\tilde{\xi} &
\end{array}.
	\]
\end{theorem}
\begin{proof}
\cite[Th. 1.6]{Art:FormModuli}
\end{proof}
\begin{cor}\label{algebraizable1}
Every deformation $\mathcal{X} \rightarrow 
 \mathrm{Spf}R$ can be extended to a deformation $\mathcal{X} \rightarrow 
\Spe R$.
\end{cor}
\begin{proof}
By lemma \ref{can-map1} the local deformation 
functor attached to a wild ramification point,
is effective, and since $\Dgl$ is the 
product of the local deformation functors
at wild ramification points \cite[3.3.4]{Be-Me}, \cite[1.10.1]{CK}
we have that $\Dgl$ is also effective.
Moreover, lemma \ref{Mpelijohn} and
proposition \ref{coherent1} imply that 
$\Dgl$ is locally of finite presentation 
and the desired result follows by Artin's
criterion of algebraisation \ref{ACOA}.
\end{proof}
%
%
%
%
%
\subsection{Construction of the $A$-free modules} \label{sub-sec3}

%
%
%
%
\begin{lemma} \label{proj-A-curves}
Let $S$ be any scheme and  let $\mathcal{X}\rightarrow S$ be
a deformation of $X$, so that the special fibre is a nonsingular curve of genus $g\geq2$. The relative curve $\mathcal{X} \rightarrow S$ is projective. 
\end{lemma}
\begin{proof}
Observe that a deformation of a nonsingular curve of genus $g\geq 2$ is a stable curve. The desired result follows by \cite[cor. p. 78]{DelMum}.
\end{proof}

\begin{lemma} \label{action-mod-lift}
Let $\mathcal{X}\rightarrow \Spe A$ be a deformation over the 
local domain $A$, let $P$ be a wild ramified point on the 
special fibre  and
        let $T=\{\bar{P}_i\}_{i=1,\cdots s}$ be the set of horizontal branch divisors that restrict 
        to $P$ in the special fibre of $\mathcal{X}$, and let $D$ be a $G$-invariant
	divisor supported on $T$ of degree $d$. There is a natural number $a$ such that the  map 
        \[
        \Gamma(\mathcal{X}, \O_{\mathcal{X}} (aD) )\otimes _A k \rightarrow L(da P)
        \]
        is an isomorphism. Then  $ \Gamma(\mathcal{X}, \O_{\mathcal{X}} (sD) )$
	is a free $A$-module, and there is a free $G$-invariant $A$-module $L$, such that 
	$L\otimes _A k= L(dP).$
\end{lemma}
\begin{proof}
The divisor $D$ is an effective Cartier divisor, and  flatness of $D$ 
over $\Spe A$ implies
\[
I(D) \otimes _A k \cong I (D \otimes _A k),
\]
{\em i.e.}, the ideal sheaf of $D$ restricted to the special fibre coincides 
with the ideal sheaf of the restriction of the divisor $D$ to the special 
fibre \cite[p. 7]{KaMa}. 

By inverting the ideal sheaves above we obtain
\[
I(D)^{-1}\otimes _A k \cong I(D\otimes _A k)^{-1} \Rightarrow \mathcal{L}(D) 
\otimes _A k \cong \mathcal{L}(D\otimes _A k).
\]

We would like to take global sections of the above two sheaves, 
in order to prove that 
\[
\Gamma(X_s, \mathcal{L}(D)) \otimes _A k = \Gamma(X_s,\mathcal{L}(D\otimes k))=L(dP).
\]
For all $i\geq 0$ there is a natural map \cite[prop. III 12.5]{Hartshorne:77}
\[
\phi_i: H^i(\mathcal{X},\mathcal{L}(D)) \otimes _A k \rightarrow 
H^i(X_s,\mathcal{L}(D\otimes k)).\]
 We are interested in global sections {\em i.e.},  for the zero cohomology groups, but in general $\phi_0$ 
 can fail to be an isomorphism. Instead of looking at $D$ we will consider
 $aD$, where $a$ is a sufficiently large natural number. The degree 
 of the divisor $aD$ remains the same in the special and in the generic 
 fibre, since  $H^0(\mathcal{X},\mathcal{L}(aD)/\mathcal{O}_\mathcal{X})$
 is a  free $A$-module of rank $\deg(aD)$, \cite[1.2.5]{KaMa}. 
 Let $\mathcal{X}_s$ and $\mathcal{X}_\eta$ denote the special and the 
 generic fibre of $\mathcal{X}$ and let $K$ be the field of quotients 
 of $A$.
 We will employ the Riemann-Roch theorem in both the special 
 and the generic fibre and we can choose $a$ sufficiently big 
 so that the index of speciality in both the generic and the special fibre is zero. 
 Thus, the Riemann-Roch theorem implies:
 \[
 \dim _k H^0(\mathcal{X}_s,\mathcal{L}(aD\otimes_A k))=
 \dim _K H^0 (\mathcal{X}_\eta, \mathcal{L}(aD \otimes K))=a \deg D +1-g.
\]
Let $f:\mathcal{X}\rightarrow \Spe A$ be the structure map. 
By lemma \ref{proj-A-curves} the $A$-curve is projective and 
Grauert theorem \cite[III.12.9]{Hartshorne:77}, implies that
$R^0 f_*(\mathcal{L}(aD)$ is a locally free sheaf on $\Spe A$, 
and quasicoherent by \cite[III.8.6]{Hartshorne:77}. Since $A$ is a local ring we have 
that $H^0(\mathcal{X},\mathcal{L}(aD))=H^0(\Spe A,f_*\mathcal{L}(aD)) $  is a free $A$-module. 

Consider the  $k$-subspace $L(dP) \subset L(adP)$. Since $D$ is $G$-invariant 
$L(dP)$ is also a $G$-invariant subspace of $L(adP)$. 
Let $\bar{x}_1,\ldots,\bar{x}_\ell$ be a basis of $L(dP)$ and 
let $x_1,\ldots,x_\ell \in H^0(\mathcal{X},\mathcal{L}(aD))$ that reduce to 
$\bar{x}_i$ modulo $m$.

The free submodule of $H^0(\mathcal{X},\mathcal{L}(aD))$ 
generated by $x_i$ is $G$-invariant and reduces to 
${L}(dP)$ modulo $m$. 

%
\end{proof}

%
%
%
Let $T$ be as in lemma \ref{action-mod-lift}.
Let $O(T)$ be the set of orbits of $T$ under the action of the group $G$,
on $T$. 
The $A$-module $\O_\mathcal{X}(D)$ is invariant under the action of $G$ 
if and only if, the divisor $D$ is of the form:
\begin{equation} \label{pp-oo13}
D=\sum_{C\in O(T)} n_C \sum_{P\in C} P, 
\end{equation}
{\em i.e.}, horizontal Cartier divisors  that are in the same orbit of
the action of $G$ must appear with the same weight in $D$. 

Given a space $L(iP)$ we would like to construct a $G$-invariant divisor $D$ supported 
on $T$ that in turn  will give a $G$-invariant $A$-module $\O_\mathcal{X}(D)$.
If $i=\sum_{C\in O(T)} n_C \# C$, where $n_C$ are non-negative
integers, then we can consider the $G$-invariant divisor given in (\ref{pp-oo13}).

We have proved the following:
\begin{pro} \label{semigroup-criter}
        If the semigroup $\sum_{C\in O(T)} n_C \# C$,  $n_C \in \mathbb{N}$, contains 
        the Weierstrass  semigroup of the branch point $P$ of the special fibre, then the 
        assumption of proposition  \ref{crit-matrix-lift} is satisfied and the representation 
        can be lifted. 
\end{pro}
\begin{cor}\label{def2matrix}
Let $A$ be a local Artin algebra that is dominated by
an integral local ring $R$, and suppose that 
the deformation $\mathcal{X}\rightarrow \Spe A$ can be lifted 
to a deformation with base $R$, such that the assumptions of proposition
\ref{crit-matrix-lift} are satisfied.  
The faithful representation defined in (\ref{fai-rep}), 
can be lifted to a representation:
\[
\rho_1:G_1(P) \rightarrow GL_n(A),
\]
such that $\sigma(g)$ is an upper  triangular matrix with $1$ at 
the diagonal.
\end{cor}
\begin{proof}
Suppose that  we have a lift $\rho_1$ 
of the faithful 
representation defined in (\ref{fai-rep}). 
There is a a surjective map $R\rightarrow A$, with kernel an ideal 
$I$ of $R$. The desired result follows by considering the 
representation matrices modulo the ideal $I$.  
\end{proof}
The following lemmata give us criteria, in order to prove that the assumption of 
of proposition \ref{semigroup-criter} is satisfied.

\begin{lemma} \label{semigroup-criter2}
If one orbit of $G$ acting on $T$ is a singleton, then
the semigroup
\[\sum_{C\in O(T)} n_C \# C, \;\;\; n_C \in \mathbb{N},\]
is the semigroup of 
natural numbers and the assumption of proposition \ref{semigroup-criter}
is satisfied.
\end{lemma}
\begin{lemma} 
 If $\#T \not\equiv 0 \mathrm{mod} p$ then there is at least 
one orbit that is a singleton. 
\end{lemma}
\begin{proof}
If all orbits have more than one element then all orbits must have 
cardinality  divisible by $p$, and since the set $T$ is the disjoint union 
of orbits it also must have cardinality divisible by $p$.
\end{proof}
\begin{lemma}
If $m$ is the first pole number that is not divisible by the characteristic, and $p\nmid m+1$ then 
there is an orbit that consists of only one element.
\end{lemma}
\begin{proof}
By  proposition \ref{bertin-gen} the 
Artin representation at the special fibre equals the sum of the 
Artin representations at the generic fibre. Thus if there is no singleton orbit then the sum of Artins representations at the generic fibre
is divisible by $p$ and this contradicts  corollary \ref{nodivart}.
\end{proof}

If the assumptions of proposition \ref{semigroup-criter} are not satisfied then we 
have to restrict to the ``good deformations".  We will describe them, following B. Mazur in \cite[p.289]{MazDef},
by defining  the notion 
of a deformation condition. 
\begin{orism}
Let $X$ be an algebraic curve defined over the algebraically closed
field $k$.
We define the category  $\mathcal{A}$  of deformations of curves
 over Artin local rings, whose objects are deformations
$
\mathcal{X} \rightarrow \Spe A
$
of the initial curve $X$, together with a fibrewise action of the 
group $G$ on $\mathcal{X}$, and the morphisms 
\[(\mathcal{X} \rightarrow \Spe A) \rightarrow
(\mathcal{X}' \rightarrow \Spe A') 
\]
are given by a local  algebra homomorphism 
$A \rightarrow A'$ and an $\Spe A'$-map 
$\phi:\mathcal{X} \times_{\Spe A} \Spe A' \rightarrow \mathcal{X}'$, 
making the following diagram commutative:
\[
\xymatrix{
\mathcal{X} \ar[d]  & \mathcal{X} \times_{\Spe A} \Spe A'  \ar[l]  \ar[d] \ar[r]^{\;\;\;\;\;\;\;\;\;\;\;\;\;\;\phi} &  \mathcal{X}' \ar[dl] \\
\Spe A  & \Spe A' \ar[l] & 
},
\]
and moreover $\phi$ induces the identity on the special fibre $X$.
\end{orism}
\begin{orism}
By a deformation condition on $\mathcal{A}$ we mean a 
full subcategory $\mathcal{DF}$ of $\mathcal{A}$ satisfying the following conditions:
\begin{enumerate}
\item For any morphism $(\mathcal{X},A) \rightarrow (\mathcal{X}',A')$
if $(\mathcal{X},A)\in Ob(\mathcal{DA})$ then 
$(\mathcal{X}',A') \in Ob(\mathcal{DA})$
\item Let $A,B,C$ be artinian $k$-algebras fitting into a diagram
\[
\xymatrix{ A  \ar[rd]_{\alpha} & & B \ar[ld]^{\beta} \\
 & C&}
\]
We denote by $A\times _C B$ the subring of 
$A \times B$ consisted of elements $(a,b)$ such that 
$\alpha(a)=\beta(b)$ \cite[p.270]{MazDef}. Moreover there are 
two projections $p_A,p_B$ 
making the following diagram commutative:
\[
\xymatrix{ 
 &A \times_C B  \ar[dl]_{p_A} \ar[dr]^{p_B}& \\
A  \ar[rd]_{\alpha} & & B \ar[ld]^{\beta} \\
 & C&}
\]
Consider an object $(\mathcal{X},A \times _C B)$ and let 
\[\mathcal{X}_A:=\mathcal{X}\times_{\Spe (A \times _C B)} \Spe A 
\rightarrow \Spe A\]
 and \[\mathcal{X}_B:=\mathcal{X}\times_{\Spe (A \times _C B)} \Spe B 
\rightarrow \Spe B\]
be the fibre products with repsect to the maps $p_A,p_B$.
We ask that $(\mathcal{X},A \times _C B) \in Ob(\mathcal{DA})$ if 
and only if both 
$\mathcal{X}_A,\mathcal{X}_B \in Ob(\mathcal{DA}$.

\item For any morphism $(\mathcal{X},A) \rightarrow (\mathcal{X}',A')$
if $(\mathcal{X}',A')\in Ob(\mathcal{DA})$ and $A\rightarrow A'$ is 
injective then 
$(\mathcal{X},A) \in Ob(\mathcal{DA})$
\end{enumerate}
\end{orism}
Suppose that we have  a deformation condition $\mathcal{DA}$ of $\mathcal{A}$.
We define a subfunctor $\mathcal{D}$ of the global deformation  functor 
$D$ of Bertin-M\`ezard 
\[
\mathcal{D} \subset D: \mathcal{C} \rightarrow \mathrm{Sets},
\]
where $\mathcal{D}(A)$ contains the elements of $D(A)$ that 
are objects of $\mathcal{DA}$.

\begin{pro} \label{hullcond}
The subfunctor $\mathcal{D}$ satisfies the three first Schlessinger 
criteria and has a hull $R_{\mathcal{D}}$. 
\end{pro}
\begin{proof}
The deformation conditions imply that $\mathcal{D}$ is relatively 
representable \cite[p.278]{MazDef}, and since the global deformation 
functor $D$ satisfies the first three Schlessinger criteria 
\cite[]{Sch},\cite[p. 277]{MazDef},
the functor $\mathcal{D}$ also satisfies them and therefore it has a hull.
\end{proof}

We are mainly interested in the deformation condition given in proposition
\ref{semigroup-criter}.
\begin{lemma}
The condition given in \ref{semigroup-criter} defines a deformation 
condition.
\end{lemma}
\begin{proof}
In order to prove the desired result we notice that 
is enough to prove that if  $d_i:\mathcal{X}_i \rightarrow \Spe A_i$ are  
 elements in $\mathcal{A}$ and 
$\phi: d_1 \rightarrow d_2$
then $d_1 \in \mathcal{DA}$ if and only if $d_2 \in \mathcal{DA}$.

Let $D$ be a $G$-invariant divisor of $\mathcal{X}_1$ such that 
$D \cap X=n P$. Then the divisor $\phi_*(D)$ is a $G$-invariant 
divisor on $\mathcal{X}_2$ such and $\phi_*(D) \cap X =n P$, since 
$\phi$ reduces to the identity on the special fibres.

If, on the other hand, $D$ is a $G$-invariant divisor of $\mathcal{X}_2$
such that $D \cap X =n P$ then $\phi^*(D)$ is a $G$-invariant divisor on 
$\mathcal{X}_1$ such that $\phi^*(D) \cap X =nP$, since $\phi$ reduces to the identity on the special fibre. 
\end{proof}

\section{Explicit Deformations of Matrix Representations}
In this section we will employ the construction for universal deformation 
rings for matrix representations,  explained by 
B. de Smit and H. W. Lenstra in \cite{SMLE:97}.
Let $G$ be a finite group with identity $e$. We denote by 
$k[G,n]$ the commutative $k$-algebra generated by $X_{ij}^g$ for 
$g\in G, 1\leq j \leq i \leq n$, such that 
\[
X_{ij}^e=\begin{cases} 1 & \mbox{ if } i=j \\ 0  & \mbox{ if } i \neq j \end{cases}
\]
\begin{equation} \label{ostru-123}
X_{ij}^{gh}= \sum_{l=1}^n X_{il}^g X_{lj}^h \mbox{ for } g,h \in G  \mbox{ and }
1 \leq i,j \leq n.
\end{equation}
\[
X_{ii}^g=1, \mbox{ for all } i=1,...,n, \mbox{ and for all } g\in G 
\]
and finally
\[
X_{ij}^g=0 \mbox{ for } i<j \mbox{ and for all } g\in G.
\]
Let $A$ be a $k$-algebra.
Consider the multiplicative group $L_n(A)<GL_n(A)$, of invertible  
lower  triangular matrices with entries in $A$, and $1$ in the 
diagonal.
We will focus on representations 
on $L_n(A)$. 
For every ${k}$-algebra $A$ we have a canonical 
bijection 
\[
Hom_{k-\mathrm{Alg}} (k[G,n],A)  \cong Hom (G, L_n (A)),
\]
where a $k$-algebra homomorphism $f:k[G,n]\rightarrow A$ 
corresponds to the group homomorphism $\rho_f$ that sends 
$g\in G$ to the matrix $(f(X_{ij}^g))$. 
The representation $\rho: G \rightarrow L_n(k)$ corresponds to 
a  homomorphism $k[G,n]\rightarrow k$. Its kernel is a maximal 
ideal, which we denote by $m_\rho$. We take the completion 
$R(G)$ of $k[G,n]$ at $m_\rho$. The canonical map $k[G,n]\rightarrow R(G)$, 
gives rise to a map $\rho_{R(G)}:G \rightarrow L_n(R(G))$, 
such that the diagram:
\[
\xymatrix{
G \ar[r]^{\rho_{R(G)}}  \ar[d]_{=} & L_n(R(G)) \ar[d] \\
G \ar[r]^{\rho} & L_n(k) 
}
\]
is commutative.

We consider the following  functor from the category $\mathcal{C}$ 
of local Artin $k$-algebras to the category of sets 
\[
F: A \in Ob(\mathcal{C}) \mapsto \left\{
\begin{array}{l}
\mbox{liftings of } \rho: G \rightarrow L_n(k) \\
\mbox{to } \rho_A: G \rightarrow L_n(A) 
\mbox{modulo} \\ \mbox{conjugation by an element }\\
\mbox{of }  
\ker(L_n(A)\rightarrow L_n(k))
\end{array}
\right\}
\]

The ring $R(G)$ defined above does not represent the 
deformation functor $F$, since $A$-equivalent deformations may correspond 
to different maps in $Hom(R(G),A)$.

If $n=2$, {\em i.e.}, in the case of a two dimensional 
representation, the conjugation action is trivial 
and $F$ is representable by  $R(G)$. 

In the deformation theory over fields of characteristic zero,
if the representation is irreducible,
one takes the closed subalgebra of $R(G)$ generated by 
the traces of elements in $\rho:G\rightarrow Gl_n(R)$, 
and since characters distinguish equivalence classes
of representations, this subalgebra represents the deformation 
functor.

We are working over fields of positive characteristic 
and this approach is not suitable: The  representation 
can be chosen  to be  indecomposable but not irreducible. 
The theory of Brauer characters, does also not help very much, 
in the case of equicharacteristic deformations, since 
Brauer characters do not take values in $R(G)$.

We will avoid to answer whether the functor $F$ is representable;
for our needs it is enough that there is a natural transformation 
from $Hom_{k-alg}(R(G),\cdot)\rightarrow F(\cdot)$.

Let $D:\mathcal{C} \rightarrow {\rm Sets}$ be the functor of 
J. Bertin, A. M\'ezard \cite{Be-Me}, introduced in (\ref{Bertin-Mezard-functor}).
We will define a natural transformation of functors
$\Psi: F \rightarrow D$,
and using this natural transformation we will prove that the Krull dimension 
of the hull of the deformation functor $D$ equals the dimension 
$\dim_k F(k[\epsilon])$, of the tangent space of the functor 
$F$. 

Let $S=\{t^n, n\in \mathbb{N}\}$ considered as a multiplicative system.
We choose elements  $F_i\in  A[[t]] S^{-1}$, $i=1,\ldots,m$ so that $F_i \equiv f_i \mod m_A$.
The elements $F_i$ can be written as 
\[
F_i =\frac{u_i(t)}{t^{\lambda(i)}}, 
\]
where $u_i(t)$ is a unit in $A[[t]]$ an $\lambda(i)$ is the pole order of the function $f_i$. 
Moreover we assume that $u_i(t)$ is  a polynomial in $A[t]$.

An element in $F(A)$, defines a linear  action on the free $A$-module generated by the elements $F_i$ given by 
\begin{equation} \label{linear-action-lift}
\sigma(F_i)=\sum_{\nu=0}^i \rho_{i\nu} (\sigma) F_\nu.
\end{equation}
The above action is not necessary compatible with  multiplication, 
{\em i.e.}, $\sigma(ab)$ and $\sigma(a) \sigma(b)$ need not to be equal.
We will give conditions so that the action is compatible with multiplication.
In order to do this we have to determine the value of $\sigma(t)$ and expand to 
the elements of $A[[t]]$ so that $\sigma$ respects addition and multiplication.
This  definition of $\sigma(t)$ should give the same results  on $F_i$ with 
(\ref{linear-action-lift}).

In order to do such a computation, we will need a general version of 
Hensel lemma. We will follow the notation of Bourbaki \cite[III 4.2]{BouComm}. We will say that a ring $R$ is 
linearly topologized, if there exists a fundamental system of neighbourhoods of $0$ consisting of ideals 
of $R$. If the commutative ring $R$ is linear topologized, Hausdorff and complete we will say that 
the ring satisfies Hensel's conditions. Let $(R,m)$ be an ordered pair, consisted of a ring $R$ and 
an ideal $m$ of $R$, so that $R$ satisfies Hensel's conditions and the ideal $m$ is closed and the elements 
of $m$ are topologically nilpotent, {\em i.e.}  for all $x\in m$ we have $\lim_{\nu \rightarrow \infty} x^\nu=0$.
Then the pair $(R,m)$ is said to satisfy Hensel's conditions \cite[III 4.5]{BouComm}.
We will denote by $R\{X\}$ the subring of $R[[X]]$ consisted of elements $\sum_{i=0}^\infty a_i X^i$ so 
that $\lim_{i\rightarrow \infty} a_i=0$.

We will need the following 
\begin{lemma} \label{Hensel1}
Let $R$ be a ring and $m$ an ideal of $R$ so that the ordered pair $(R,m)$ satisfies Hensel's conditions.
Let $f\in R\{X\}$, $a\in R$ end write $e=f'(a)$. If $f(a)\equiv 0 \mod e^2m$, then there exists $b\in R$ such that 
$f(b)=0$ and $b\equiv a \mod em$. If $b'$ is another element of $R$ such that $f(b')=0$ and $b'\equiv a \mod em$ 
then $e(b'-b)=0$. In particular, $b$ is unique if $e$ is not a divisor of zero in $R$.
\end{lemma}
\begin{proof}
Corollary 1 \cite[III 4.5]{BouComm}
\end{proof}

\begin{lemma}
Let $A$ be an Artin local ring with maximal ideal $m_A$. We consider the ring $A[[t]]$ so that the ideals 
$\langle m_A^\mu , t^\nu A[[t]] \rangle$ form a fundamental system of  neighbourhoods of $0$.
Then the ordered pair $(A[[t]],m_A)$ satisfies Hensel's conditions.
\end{lemma}
\begin{proof}
The ring $A[[t]]$ is by construction linear topologized. The ideal $m_A$ is an element of the the system 
of fundamental neighbourhoods of $0$ so it is open and closed. 
Moreover the ring $A[[t]]$ is Hausdorff and complete.
\end{proof}

Let us consider the   last row of (\ref{linear-action-lift}):
\[
\sigma(F_m)=\sum_{\nu=1}^m \rho_{m\nu} (\sigma) F_\nu(t).
\]
Let  $Y$ be the desired value of the the lift of $\sigma(t)$. Then the last equality can be written as 
\[
\frac{u_m(Y)}{Y^m}=\sum_{\nu=1}^m \rho_{m\nu}(\sigma) F_\nu(t) \Rightarrow 
\]
\begin{equation} \label{last-row2}
t^m u_m(Y)-Y^m \sum_{\nu=1}^m \rho_{m\nu}(\sigma) F_\nu (t) t^m=0.
\end{equation}
If we consider (\ref{last-row2}) modulo $m_A$ then it has as solution the value $\sigma(t)\in \A k[[t]]$.
We will use lemma \ref{Hensel1} with $a=\sigma(t)$. The derivative of (\ref{last-row2}) is 
computed
\[
t^m \left( \frac{\partial u_m(Y)} {\partial Y} - m Y^{m-1} \sum_{\nu=1}^m \rho_{m\nu}(\sigma) F_\nu (t) \right),
\]
and 
\[
e=t^m \left( \frac{\partial u_m(\sigma(t))} {\partial Y} - m \sigma(t)^{m-1} \sum_{\nu=1}^m \rho_{m\nu}(\sigma) F_\nu (t) \right).
\]
Observe that the zero divisors in $A[[t]]$ have coefficients in the maximal ideal  $m_A$, therefore
since  $e \mod m_A$ is computed to be $-m t^m \sigma(t)^{m-2}$ and 
it is not a zero divisor, $e$ is not a zero divisor in $A[[t]]$ as well. Thus, lemma \ref{Hensel1} implies that 
there is a unique solution $b$ of (\ref{last-row2}) that reduces to $\sigma(t)$ modulo $m_A$.
The value $b$ is a candidate for the lift $\tilde{\sigma}(t)$.
We observe that the equations (\ref{linear-action-lift}) give also information for the value of the 
lift of $\sigma(t)$.

We say that the tuple $(F_1,\ldots,F_m, \rho:G \rightarrow GL_n(A) )$ is 
compatible if the unique solution $b$ of (\ref{last-row2}) satisfies the  equations (\ref{linear-action-lift})  for $1\leq i <m$. 
Of course we know that $(f_1,\ldots,f_m,\rho:G \rightarrow GL_n(k))$ is compatible.
\begin{lemma} \label{equiv-compat}
If the tuple $(F_1,\ldots,F_m, \rho:G \rightarrow GL_n(A) )$ is compatible and 
$Q\in L_n(A)$ then the tuple 
$(Q F_1,\ldots,Q F_m, \rho:G \rightarrow Q GL_n(A) Q^{-1} )$
is also compatible.
\end{lemma}
\begin{proof}
For every $1\leq i \leq m$ we have $Q F_i =\sum_{\mu=1}^m Q_{i \mu} F_\mu.$
The system of equations (\ref{linear-action-lift}) can be written in matrix form as
\[
\begin{pmatrix}
F_1(Y) \\ \vdots \\ F_m(Y) 
\end{pmatrix}
=
\rho(\sigma) 
\begin{pmatrix}
F_1(t) \\ \vdots \\ F_m(t) 
\end{pmatrix}.
\]
Thus, 
\[
Q \cdot 
\begin{pmatrix}
F_1(Y) \\ \vdots \\ F_m(Y) 
\end{pmatrix}
=
Q \rho(\sigma) Q^{-1} \cdot Q 
\begin{pmatrix}
F_1(t) \\ \vdots \\ F_m(t) 
\end{pmatrix},
\]
and the desired result follows.
\end{proof}

\begin{lemma}
Let $A$ be an artin local ring with maximal ideal $m_A$ so that $A/m_A=k$.
Let $(F_1,\ldots,F_m,\rho:G \rightarrow L_n(A))$ be a compatible tuple, and for every 
$\sigma \in G$ let $\tilde{\sigma}$ denote the corresponding automorphism in $\A A[[t]]$.
For every $\sigma_1,\sigma_2\in G$ we have 
\[
\tilde{\sigma}_1\tilde{\sigma}_2(t)=\widetilde{\sigma_2\sigma_1}(t).
\] 
 \end{lemma}
\begin{proof}
The element $b_{\sigma_2}=\tilde{\sigma}_2$ is a root of 
\[
t^m u_m(Y)-t^m Y^m \sum_{\nu=1}^m \rho_{m\nu} ({\sigma}_2) F_\nu(t).
\] 
By applying $\tilde{\sigma}_1$ to the above equation we obtain
\[0=
\tilde{\sigma}_1(t) \left( 
u_m (\tilde{\sigma}_1(b_{\sigma_2})) -\tilde{\sigma}_1(b_{\sigma_2})^m  \sum_{\nu=1}^m \rho_{m\nu}(\sigma_2) \tilde{\sigma}_1 F_\nu(t)
\right)=
\]
\[
=\tilde{\sigma}_1(t) \left( 
u_m (\tilde{\sigma}_1(b_{\sigma_2})) -\tilde{\sigma}_1(b_{\sigma_2})^m  
\sum_{\nu=1}^m \rho_{m\nu}(\sigma_2) \sum_{\mu=1}^\nu \rho_{\nu\mu}(\sigma_1) F_\mu(t)=
\right)=
\]
\[
=\tilde{\sigma}_1(t) \left( 
u_m (\tilde{\sigma}_1(b_{\sigma_2})) -\tilde{\sigma}_1(b_{\sigma_2})^m  
\sum_{\nu=1}^m \rho_{mu\nu}(\sigma_2 \sigma_1) F_\mu(t) \right).
\]
Since the element $\tilde{\sigma}_1(t)$ is not a zero divisor we obtain that
\[
\tilde{\sigma}_1 (\tilde{\sigma}_2(t))
\]
is the unique root of 
\[
t^m u_m(Y)-t^m Y^m \sum_{\nu=1}^m \rho_{m\nu} ({\sigma_2 \sigma_1}) F_\nu(t).
\] 
and the desired result follows.
\end{proof}

\begin{orism}
We will say that the tuples $(F_1,\ldots,F_m, \rho:G \rightarrow GL_n(A) )$ and 
$(F_1',\ldots,F'_m, \rho':G \rightarrow GL_n(A) )$ are equivalent if there is an 
element $Q\in L_n(A)$ so that 
\[(F_1',\ldots,F'_m, \rho':G \rightarrow GL_n(A) )=(Q F_1,\ldots,Q F_m, \rho:G \rightarrow Q GL_n(A) Q^{-1} )\]
We define the functors
\[
\mathcal{F}_1: A \in Ob(\mathcal{C}) \mapsto \left\{
\begin{array}{l}
\mbox{Set of   tuples  } \\
\big(F_1,\ldots,F_m, \rho:G \rightarrow GL_n(A) \big)\\
\mbox{ defined  over } A
\end{array}
\right\}
\] 
and
\[
\mathcal{F}: A \in Ob(\mathcal{C}) \mapsto \left\{
\begin{array}{l}
\mbox{equivalence classes } \\
\mbox{of tuples over } A
\end{array}
\right\}
\]
If $g:B\rightarrow A$ is a morphism of local artin algebras and 
\[
\Omega:=(F_1,\ldots,F_m,\rho:G\rightarrow L_n(B)) \in \mathcal{F}_1(B)
\]
\[\mathcal{F}_1(g)(\Omega)=
(g(F_1),\ldots,g(F_m),g\circ\rho:G\rightarrow L_n(A)) \in \mathcal{F}_1(A).
\]
where $g\circ \rho(\sigma)$ denotes the matrix $g(\rho_{\mu\nu}(\sigma))$, {\em i.e.},
$g$ is applied at every entry of the matrix $\rho(\sigma)$. Similarly, we can define the map
$\mathcal{F}(g): \mathcal{F}(B) \rightarrow \mathcal{F}(A)$.
\end{orism}
Observe that compatibility of equivalence classes of tuples is well defined  by lemma  \ref{equiv-compat}.
\begin{pro}
Compatibility of tuples is a deformation condition in both the functors $\mathcal{F},\mathcal{F}_1$.
\end{pro}
\begin{proof}
Let us consider a pair $(A,\Omega)$, consisted of an artin local 
ring $A$ together with a tuple $ \Omega:=(F_1,\ldots,F_m, \rho:G \rightarrow L_n(A))$. Let us consider the 
set of compatible tuples over $A$ by $\mathcal{F}_1'(A) \subset \mathcal{F}_1(A)$.
 A morphism  $\phi:A \rightarrow A'$ of local artin algebras, induces  the map $\Omega \rightarrow \mathcal{F}_1(\phi)(\Omega)$, 
on tuples.
In order to have a deformation condition we should check that 
\begin{enumerate}
\item 
For every local artin algebras $A,A'$ and morphism $f:A \rightarrow A'$, if $\Omega \in \mathcal{F}_1'(A)$ then 
$\mathcal{F}_1(\phi)(\Omega)\in \mathcal{F}_1'(A')$.
\item 
Let  $A,B,C$ be  local artin algebras and $\alpha:A \rightarrow C$,  $\beta:B \rightarrow C$ be  two morphisms
of local artin algebras. We form the local artin algebra 
\[A\times _C B=\{ (x,y): x\in A, y\in B, \alpha(x)=\beta(y)\},
\]
equipped with coordinate addition, multiplication, and scalar $k$-multiplication.
Let $\Omega \in \mathcal{F}'_1 (A\times _C B)$ we want the following
\[
\Omega \in \mathcal{F}_1' \Leftrightarrow \mathcal{F}_1(\alpha)(\Omega)  \in \mathcal{F}_1'(A) 
\mbox{ and } \mathcal{F}_1(\beta)(\Omega) \in \mathcal{F}_1'(B)
\]
\item 
For every local artin algebras $A,A'$ and every injective map $f:A \rightarrow A'$, we want the following
\[
\Omega \in \mathcal{F}_1'(A) \Leftrightarrow \mathcal{F}_1(\phi)(\Omega)\in \mathcal{F}_1'(A').\]
\end{enumerate}
Checking all this conditions is straightforward. A similar check can be done also for the functor $\mathcal{F}$.
\end{proof}

Consider the functor $Hom_{k-alg} (R(G),\cdot)$. Every fixed selection $(F_1,\ldots,F_m)$,
$F_i\in A[[t]] S^{-1}$, $S:=\{t^n, n \in \mathbb{N} \}$  gives rise to 
a natural map of $Hom_{a-alg} (R(G),\cdot) \rightarrow \mathcal{F}_1$ and $F(\cdot) \rightarrow \mathcal{F}$.

The compatability of tuples gives rise to a set of  deformation conditions (that depend on the basis 
selection $F_1,\ldots,F_m)$ on $Hom_{k-alg}(R(G),\cdot)$. By lemma on page 279 of \cite{MazDef}, we have that 
the deformation conditions on $Hom_{k-alg}(R(G),\cdot)$ give rise to a representable functor
$Hom(R(G)/I,\cdot)$, represented by the $k$-algebra $R(G)/I$, where $I$ is a suitable ideal, 
that depends on the selection of $F_1,\ldots,F_m$.

\begin{lemma}
The ideal $I$ is an ideal that is generated by polynomials on the variables $X_{i\nu}^g$, $g\in G$, $1\leq i <m$, 
{\em i.e.} the deformation condition does not involve the variables $X_{m,\nu}^g$.
\end{lemma}
\begin{proof}
The set of compatibility conditions (\ref{linear-action-lift}) give the following conditions on $R(G)$:
For every $g$ the unique solution $Y$ of 
\[
u_m(Y) =Y^m \sum_{\nu=1}^m X_{i,\nu}^g X_{i \nu}^g 
\]
should satisfy the equations 
\[
u_i(Y) =Y^i \sum_{\nu=1}^i X_{i\nu}^g F_{\nu}(t),
\]
 for every $1\leq i <m$.
\end{proof}

\begin{pro} \label{relateFD}
There is a natural transformation  $\Psi$
from the functor $\mathcal{F}'$ to the functor ${D}$.
\end{pro}
\begin{proof}
Let $A \in Ob(\mathcal{C})$ be a local Artin algebra with 
maximal ideal $m_A$.
Consider the $n$ dimensional $k$-vector space $L(mP)$, 
together with the flag $L(iP), i< m$ of vector spaces. 

Every tuple in $\mathcal{F}(A)$ gives rise to a sequence $L_i$ of free $A$-modules
where 
\[{L_i:=_A \langle F_1,\ldots, F_i \rangle},\]
 so that $L_i \otimes_ A k =L(iP)$  and an action of $G$ on them
that is a lift of the action of $G$ on the spaces $L(iP)$.

Since the tuple $\Omega:=(F_1,\ldots,F_m, \rho:G \rightarrow L_n(A) )$ is compatible 
we get an automorphism $\tilde{\sigma}(t)$ that defines an equivalence class  $\Phi(\Omega)\in D(A)$.

Let $A,B$ we local artin algebras and consider a morphism   $g:B\rightarrow A$. For every 
tuple $ \Omega:=(F_1,\ldots,F_m, \rho:G \rightarrow L_n(B)) \in \mathcal{F}(B)  $ the map $\mathcal{F}(g)$ gives the tuple
$\mathcal{F}(g) (\Omega) \in \mathcal{F}(A)$.

If $\tilde{\sigma}$ denotes the automorphism of $\A B[[t]]$ that corresponds to $\Omega$ then 
$g(\tilde{\sigma}(t))$ is the automorphism that corresponds to $\mathcal{F}(g)(\Omega)$.

\end{proof}
{\bf Example} Consider a  functor that  is represented by a 
complete ring that is not a domain. For example let $R=k[[x_1,x_2]]/x_1^3x_2^4$. 
Then the tangent space $\mathrm{Hom}(R,k[\epsilon]/\epsilon^2)$ is two dimensional 
generated by the maps $\theta_i(x_j)=\delta_{ij} \epsilon$, but the arbitrary 
linear combination $\theta=\lambda_1 \theta_1 + \lambda_2 \theta_2$ 
 could not lift to a homomorphism $\bar{\theta}: R \rightarrow k[[\epsilon]]$,
since $\theta (x_1)^3 \theta( x_2)^4=0$ and $k[[\epsilon]]$ is a domain.

In case $R$ is a domain the following lemma shows that infinitesimal 
deformations are unobstructed:
\begin{lemma} \label{extensintdom}
Let $R$ be  a complete local $k$-domain, and denote by $m$ the 
maximal ideal of $R$. Suppose that $k=R/m$ is algebraically closed.
  Let $A$ be an Artin local ring, an suppose that 
$A=k[[a_1,\ldots,a_s]]/I$. Every 
element in $\mathrm{Hom}_k(R,A)$ can be lifted to $\mathrm{Hom}_k(R,k[[a_1,\ldots,a_s]])$.
\end{lemma}
\begin{proof}
By Noether's normalisation lemma \cite[cor. 16.18]{Eisenbud:95} there are elements
$x_1,\ldots,x_d \in R$ such that $R/k[[x_1,\ldots,x_d]]$
is a separating integral extension. 
A function $g:k[[x_1,\ldots,x_d]]\rightarrow A$, given 
by $g(x_i)=f_i(a_1,\ldots,a_s)$ can obviously extend to a 
function $g:k[[x_1,\ldots,x_d]]\rightarrow k[[a_1,\ldots,a_s]]$.

Now if $y\in R $, is an arbitrary element satisfying a separable  equation $\sum_{i=0}^n b_i y^i$, 
$b_i\in k[[x_1,\ldots,x_n]]$, then Hensel's lemma implies that the equation 
$\sum_{i=0}^n g(b_i) T^i$ has a unique solution $T$, extending the solution 
of $\sum_{i=0}^n \overline{g(b)}_i T^i$ of the equation taken modulo $m_A$; 
recall that the field $k$ is algebraically closed, and $g(t)=t$ for all $t \in k$. 
We define this unique solution to be the value of $g$ at $y$. 
\end{proof}

Using the criteria of Schlesinger \cite{Sch},\cite[p.277]{MazDef} J.Bertin A. M\'ezard  
\cite{Be-Me} were able to prove that the deformation $D$ admits a hull.  We have seen in \ref{hullcond} 
that the subfunctor $\mathcal{D}$ also admits a hull.

This means that there is a complete ring $\tilde{R}$ and a smooth 
map of functors 
$\mathcal{D}(\cdot)\rightarrow  Hom_{k-alg} (\tilde{R},\cdot) $, 
that induces an isomorphism on the tangent spaces.
The hull $\tilde{R}$ might not be a domain. 
In order to study its dimension we factor out nillpotent elements obtaining 
 $\tilde{R}/\sqrt{0}$, a ring that corresponds to
a variety. We will study the  rings of the irreducible 
components. We observe first that the rings $\tilde{R}$ and $\tilde{R}/\sqrt{0}$ 
have the same Krull dimension, and hence the Krull dimension of 
$\tilde{R}$ equals the maximum dimension of the rings that correspond to
the irreducible components of $\tilde{R}/\sqrt{0}$.

Let $R_i$ be   a ring, corresponding to an irreducible component of $\tilde{R}/\sqrt{0}$. There is an onto map 
$\tilde{R} \rightarrow R_i$, that gives rise to an injection
\[\mathrm{Hom}_{k-alg}(R_i,\cdot) \rightarrow \mathrm{Hom}_{k-alg}(\tilde{R},\cdot).\] 
We will use that fact that $R_i$ is a domain together with 
 proposition \ref{crit-matrix-lift} in order to construct a 
 map $\mathrm{Hom}_{k-alg}(R_i,\cdot) \rightarrow F(\cdot)$.
We will need the following:
\begin{lemma} \label{lemseqext12}
Suppose that $R=R_i$ is one of the rings defined above.
Let $A$ be an Artin local $k$-algebra and suppose that 
$A=k[[a_1,\ldots,a_s]]/I$. There is a  noetherian complete local domain $k[[a_1,\ldots,a_s]]/I' $, a  sequence of small extensions
\begin{equation} \label{ext:seq12}
k[\epsilon ]/\langle \epsilon ^2 \rangle  =A_0 
\stackrel{\phi_1}{\longleftarrow} A_1 \leftarrow  \cdots 
\stackrel{\phi_i}{\longleftarrow}  A_i \stackrel{\phi_{i+1}}{\longleftarrow} A_{i+1} \leftarrow  \cdots
\leftarrow
k[[a_1,\ldots,a_s]]/I',
\end{equation}
such that $A=A_n$ for some natural $n$,
and the following diagram is 
commutative:
\[
\xymatrix{
\mathcal{D}(A_{i+1}) \ar[d]_{\mathcal{D}(\phi_{i+1})} \ar[r] & \mathrm{Hom}(R,A_{i+1}) \ar[d]^{\mathrm{Hom}(\phi_{i+1})} \\
\mathcal{D}(A_i)        \ar[r]              &  \mathrm{Hom}(R,A_i)
}
\]
Moreover, for every element  $\gamma \in \mathrm{Hom}(R,A)$ there are elements $\gamma_i$  in $\mathrm{Hom}(R,A_i)$ 
and $\delta_i \in \mathcal{D}(A_i)$  such that 
$\mathrm{Hom}(\phi_{i+1})(\gamma_{i+1})=\gamma_i$ and 
$\mathcal{D}(\phi_{i+1})(\delta_{i+1})=\delta_i$.
 
\end{lemma}
\begin{proof}
Let $A$ be  an Artin local ring expressed as $A=k[[a_1,\ldots,a_s]]/I$.
The existence of the  sequence (\ref{ext:seq12}) is clear. Lemma \ref{extensintdom}
implies the existence of the sequence of elements  $\gamma_i$ that are compatible with the maps
$\mathrm{Hom}(\phi_i)$. The existence of the elements $\delta_i$ compatible with the 
elements $\mathcal{D}(\phi)$ is implied by smoothness of the map 
$\mathcal{D}(\cdot) \rightarrow \mathrm{Hom}(R,\cdot)$.
\end{proof}

\begin{lemma}
There is a natural transformation $q:\mathrm{Hom}_{k-alg}(R_i,\cdot) \rightarrow \mathcal{F}'(\cdot)$.
\end{lemma}
\begin{proof}
Start with an element $\gamma$ in $\mathrm{Hom}(R_i,A)$ for an Artin local ring $A$. Construct a sequence 
of small extensions and maps as in equation  (\ref{ext:seq12}) of lemma \ref{lemseqext12}.
By lemma \ref{lemseqext12}  there is a sequence of deformations $\delta_i \in \mathcal{D}(A_i)$ 
that leads to a deformation in $\hat{\mathcal{D}}(k[[a_1,\ldots,a_s]]/I')$, where $k[[a_1,\ldots,a_s]]/I')$ is a
notherian local complete domain. By corollary 
\ref{algebraizable1} we obtain a deformation $\mathcal{X}\rightarrow \Spe k[[a_1,\ldots,a_s]]/I'$, \
and using corollary  \ref{def2matrix} we obtain an element  $q(\gamma) \in \mathcal{F}'(A)$.
\end{proof}
We have the following picture:
\[
\xymatrix{
    &	\mathrm{Hom}_{k-alg}(R_i,\cdot) 	\ar[d]^q \ar@{^{(}->}[r]  &    \mathrm{Hom}_{k-alg}(\tilde{R},\cdot) \\
  \mathcal{F}(\cdot) &  \mathcal{F}'(\cdot) \ar@{_{(}->}[l] \ar[r] & \mathcal{D}(\cdot) 
\ar[u]_{\mbox{smooth}}
}
\]
The above  diagram induces the following diagram on the tangent spaces:
\begin{equation}\label{diag11}
\xymatrix{
   			\mathrm{Hom}_{k-alg}(R_i,k[\epsilon]) 	\ar[d]_{q} \ar@{^{(}->}[r]^a  
\ar@{^{(}->}[dr]^{\ell} 
&    \mathrm{Hom}_{k-alg}(\tilde{R},k[\epsilon])   \\
 \mathcal{F}'(k[\epsilon]) \ar[r]^{\psi} & \mathcal{D}(k[\epsilon]) \ar[u]^{\cong}_{b}
}
\end{equation}
Where $l=b^{-1}\circ a$ is an injection. 
\begin{pro}\label{prop123}
	The map $\psi$, induces an isomorphism 
	between the $k$ vector spaces $\mathcal{F}'(k[\epsilon])$ and $\mathrm{Im}(\ell)$. 
\end{pro}
\begin{proof}
	By construction $\ell$ induces an isomorphism between 
	$\mathrm{Hom}_{k-alg} (R_i,k[\epsilon])$ and $\mathrm{Im}(\ell)$. 
	Let $v \in \mathrm{Im}(\ell)$, we find the  inverse $\ell^{-1}(v)$
	and take $q(\ell^{-1}(v))$. Since the diagram (\ref{diag11}) is 
	commutative, the map $\psi$ is onto $ \mathrm{Im}(\ell)$. 

	We have seen that Hensel's lemma implies that every deformation 
	in $\mathcal{F}'(A)$ gives rise to a unique deformation in $\mathcal{D}(A)$. 
	This proves that the map $\psi$ is 1-1 and moreover 
	$\mathrm{Im}(\psi)=\mathrm{Im}(\ell)$.
\end{proof}

\begin{cor} \label{cor123}
	If the ring $\tilde{R}/\sqrt{0}$ is not a domain, then all 
	rings $R_i$ that correspond to irreducible components of
	$\Spe R$, have the same dimension, equal to 
	$\dim_k \mathcal{F}'(k[\epsilon])=\dim \tilde{R}$. 
	Moreover every infinitesimal deformation in $\psi (\mathcal{F}'(k[\epsilon]))$ 
	is unobstructed. 
\end{cor}
\begin{proof}
	We have shown that $\dim_k \mathcal{F}'(k[\epsilon])=\dim_k \mathrm{Im}(\ell)$, 
	and $\dim_k \mathcal{F}'(k[\epsilon])$ does not depend on $R_i$. On the 
	other hand every infinitesimal deformation in $\psi (\mathcal{F}'(k[\epsilon]))$, 
	is the image of the representable functor $\mathrm{Hom}_{k-alg}(R_i,\cdot)$, 
	and this implies the desired result. 
\end{proof}
	
%
%
%
%

{\bf Examples:}
{\bf 1.} $\mathbf{n=2}$. This is always the case if $p>2g-2$, because the first pole 
number is smaller than $2g-2$. Automorphisms of curves with this 
property were studied by P. Roquette in \cite{Roq:70}. 

The group  $G_1(P)$ is an elementary Abelian group 
isomorphic to a subgroup of $\{ (a_{ij})_{i,j=1,2}: a_{11}=a_{22}=1, a_{12}=0\}.$

We will check whether the assumptions of the lifting criterion of 
corollary \ref{semigroup-criter2} are satisfied. 
If the automorphism group $G$ is cyclic of order $p$, then every irreducible component 
of the ramification divisor is by construction $G$-invariant, and 
the orbits are singletons.  
If the group $G$ is elementary Abelian with more than one
cyclic components of order $p$, 
and there are no singleton orbits, then all orbits have order 
divisible by $p$. By comparing the Artin representations 
at the special and generic fibres using proposition \ref{bertin-gen}, 
we obtain that the order function of the ramification filtration is 
a multiple of $p$. 
Thus,  $p | m+1$, a contradiction. Thus the 
functors $\mathcal{D},D$ are equal.

Since the representation is two dimensional the  conjugation 
action on the ring $R(G)$ is trivial, and the functor $F$ is 
representable by 
\[
{R(G)}:= k[[ t_1,\cdots,t_{\log_p(|G_1(P)|)}  ]].
\]
The dimension of  the tangent space of $F$ and the krull dimension 
of the hull is $\log_p (|G_1(P)|)$.

\textbf{Remark:}
A comparison of this result with the computation of J. Bertin, A.M\`ezard 
for deformations of the cyclic group $\mathbb{Z}_p$
\cite[Th. 5.3.3]{Be-Me} gives us the nontrivial result: {\em  for  the smallest 
$m$ such  that $\dim L(mP)=2$, $m<p-1$}.  

%
%
%
%
%
{\bf 2.} $\mathbf{n=3}$,{\bf Ordinary Curves.} 
We will use the tools developed so far in order to study 
ordinary curves with  $n=3$ and compare our  results
with the more general  results of G. Cornelissen- F. Kato \cite{CK} (where the 
decomposition groups could have elements prime to $p$ as well).


We will use the notation from example 3. of page \pageref{OrdCurves}.
Let $\rho_{ij}(g)$ be the representation matrix. We have $\rho_{32}(g)= c_1(g)$ 
and $\rho_{21}(g)=\lambda c_1(g)$.

We form the ring $R(G)$, generated by $X_{21}^{g_i},X_{31}^{g_i},X_{32}^{g_i}$. 
We observe that there are the relations $X_{21}^{g_i}=\Lambda X_{31}^{g_i}$, 
$\Lambda \in R(G), \Lambda \equiv \lambda \mathrm{mod} m_{R(G)}$. 
Since $X_{32}^g=0 \Rightarrow X_{31}^g=0$
$X_{31}^{g_i} \in \langle X_{32}^{g_i} \rangle$,
{\em i.e.}, 
\begin{equation} \label{x31equ}
 X_{31}^{g_i} = X_{32}^{g_i} 
f(X_{ij}^{g_\nu}).
\end{equation}
Thus the only independent elements are $X_{32}^{g_i}$. 

We will now study the image of
$\mathrm{Hom}_{k-alg}(R(G),k[\epsilon])\rightarrow F(k[\epsilon])$, 
{\em i.e.}, we will study representations to $L_n(k[\epsilon])$
modulo the conjugation action. 
Let 
\[
\bar{\rho}^i(g): G_1(P) \rightarrow L_n(k[\epsilon])\;\; i=1,2
\]
be two equivalent  representations, with representation matrices
\[
\begin{pmatrix}
	1 & 0 & 0 \\
	\rho_{21}(g) & 1 & 0 \\
	\rho_{31}(g) & \rho_{32}(g) & 1 \\
\end{pmatrix}
+ 
\epsilon 
\begin{pmatrix}
	0 & 0 & 0 \\
	a_{21}^{i}(g) & 0 & 0 \\
	a_{31}^{i}(g) & a_{32}^{i}(g) & 0 \\
\end{pmatrix}.
\]
Let $Q=\begin{pmatrix}
	1 & 0 & 0 \\
	a\epsilon & 1 & 0 \\
	c\epsilon & \epsilon b & 1 \\
\end{pmatrix}$ be a matrix that transforms by conjugation 
$\bar{\rho}^1$ to $\bar{\rho}^2$.

This is equivalent (recall that $\epsilon^2=0$) to 
\begin{equation} \label{eq1diag}
	a_{21}^1(g)=a_{21}^2(g),\;\; a_{32}^1(g)=a_{32}^2(g),
\end{equation}
and 
\begin{equation} \label{eq2diag}
	a_{31}^1(g) +  b \rho_{21}^1(g) -a \rho_{32}^1(g)= a_{31}^2(g).
\end{equation}
We express equation  (\ref{x31equ}) as 
\begin{equation} \label{lala32}
X_{31}^{g_i}=X_{32}^{g_i}
\left(
a_0 + a_1 X_{31}^{g_i} +a_2 {X_{31}^{g_i}}^2+\cdots
\right),
\end{equation}
where $a_i$ are polynomials in all variables except $X_{31}^{g_i}$.
The representations $\bar{\rho}^i(g)$ are images of maps 
$f_i\in \mathrm{Hom}_{k-alg}(R(G),k[\epsilon])$. By taking the functions $f_i$ on (\ref{lala32}) we compute (recall that $\epsilon^2=0$) :
\[
a_{31}^i(g)
\left(
1-\rho_{32}(g) 
\sum_{i=1}^{deg(f)} a_i  i \rho_{31}(g)^{i-1} 
\right)= a_{32}^1(g) \sum_{i=0}^{\deg(f)} a_i \rho_{31}(g_i)^i. 
\] 
Therefore, if $\frac{\partial f}{\partial X_{31}^{g_i}}\big(\rho_{31}(g_i)\big)\neq 1$, then  
the value $a_{31}^i$ can be computed from the values of 
$a_{21}^{g_i}$ and $\rho_{ij}(g_i)$. In particular, $a_{31}^1(g)=a_{31}^2(g)$. 
If on the other hand $\frac{\partial f}{\partial X_{31}^{g_i}}\big(\rho_{31}(g_i)\big)= 1$, then 
the value $a_{31}^{g_i}$ is independ of the values 
$a_{21}^{g_i}$ and $\rho_{ij}(g_i)$. This, is not 
possible since the map $c_1$ defined in (\ref{c1})
considered  on the generic fibre 
has to be faithfull.

We identify the set of representations $G_1(P) \rightarrow L_n(k[\epsilon])$, 
with the set of homomorphisms of additive groups 
$\mathrm{Hom}( G_1(P),k)$ which is of dimension $\log_p | G_1(P)|$.
By (\ref{eq2diag}) we have  
$ 0=b \rho_{21}(g) -a \rho_{32}(g)=(b\lambda-a) \rho_{32}(g)$, and since $a,b$ are 
arbitrary this means that $\rho_{32}(g)=0$ and the  tangent space of $F(k[\epsilon])$ is identified by 
$\mathrm{Hom} (G_1(P),k) / \rho_{32}(g)$, {\em i.e.} the 
group homomorphism $\rho_{32}(g): G_1(P) \rightarrow k$ is 
considered to be zero. This is a space of dimension 
$\log_p | G_1(P)|-1$, and this result coincides with the
result of G. Cornelissen, F. Kato. 

{\bf 3.} {\bf $\mathbf{p}$-cyclic covers of the affine line.}
We consider in this case  deformations of the curve 
defined in example 4. of page \pageref{p-cycl}.
So we consider the curve $X:w^p-w=x^{p^s+1}$, and let $G$ denote the $p$-part of the 
automorphism group of $X$.
In this section we will consider deformations 
$\mathcal{X}\rightarrow \Spe A$  
of the couple  $(X,G)$, where $A$ is a local integral domain.

Consider the Galois group 
 $H':=\mathrm{Gal}(X/\mathbb{P}^1_k)=\mathbb{Z}/p\mathbb{Z}$, a cyclic group of order $p$.
Denote by  $Def_G(X), Def_H'(X)$ the global deformation functors of the curve $X$ with respect to the 
groups $G,H'$ respectively. The group $G$ acts freely on the completment of the generic fibre of $\mathcal{X}\rightarrow \Spe A$,
and  J. Bertin, A. M\'ezard in \cite{BeMe2004} proved that there is a well defined map  $\mathrm{ind}:\mathrm{Def}_G(X) \rightarrow \mathrm{Def}_{G/H'}(X/H')$.
That means that  $\mathrm{ind}\big(\mathcal{X}/H' \rightarrow \Spe A\big)\in\mathrm{Def}_{G/H'}(X/H') $ is a deformation of $X/H'=\mathbb{P}^1$. But such a deformation
is trivial {\em i.e.} $\mathcal{X}/H'=\mathbb{P}^1_A \rightarrow A$. Therefore, the generic fibre $\mathcal{X}_\eta \rightarrow \Spe Quo(A)=\Spe k$, 
is a cyclic extension of $\mathbb{P}^1$, and is given in terms of an Artin-Schreier extension:
there are 
functions $W,X$  on the generic fibre of $\mathcal{X}$, 
such that 
\[W^p -W=f({X}),\]
and $f(X) \in k(X).$
Let us write $f(x)={\prod f_i}/{\prod g_i}$, where $f_i,g_i$ are irreducible polynomials in $k[x]$ and since $k$ is assumed 
algebraically closed all are of degree $\leq 1$. 
The places ramified in the cover $\mathcal{X}_\eta \rightarrow \mathbb{P}^1$ correspond to the poles of the function $f$.
If there  are polynomials $g_i$ of degree $1$,  then there should be polynomials $f_i$ so that $g_i=f_i \mod m_A $, 
since at the special fibre only $\infty$ is ramified. But this situation does not give as a proper relative curve, since any point on the 
generic fibre gives rise to a unique horizontal divisor intersecting the special fibre at a unique point.

Thus, by  deforming the curve $X$ we deform the polynomial $f_0(x)=x^{p^s+1}$, 
and since the genus of every fibre, that is depended on the degree of 
$f$, has to remain constant, we deform $f$ by adding lower degree terms
$a_iX^i$, $i< p^s+1$, such that $a_i=0 \mod m $. 

Let $G_{1}(\infty)$ be the $p$-part of the decomposition group at $\infty$. 
According to  \cite[prop. 8.1]{CL-MM} we have  $Ad_{x^{p^s+1}}(Y)=Y+Y^{p^{2s}}$, and 
$Ad_f(Y)=\sum_{0 \leq \lambda \leq s} t_{i}^{p^{s-i}} Y^{p^{s-i}}+ t_i^{p^s} Y^{p^{s+i}}$ thus the  
group $H:=G_1(\infty)/\mathrm{Gal}(X/\mathbb{P}^1_k)$ is an elementary abelian group of order $p^{2s}$.
 The group $G_1(\infty)$ is an extraspecial group, 
{\em i.e.}, an extension of a $p$-cyclic group by an elementary abelian group. For the center we have
$Z(G)=\mathrm{Gal}(X/\mathbb{P}^1_k)=\mathbb{Z}/p\mathbb{Z}$.
Let us fix elements  $\mathcal{B}=\{g_i, i=1,\ldots,2s\}$ in $G_1(\infty)$ such that their image 
in  $H$ form a basis of $H$ 
as an $\mathbb{Z}/p\mathbb{Z}$-module of rank $2s$.
The elements $g_i$ generate $G_1(\infty)$. Indeed, the group $\langle g_i \rangle$ 
generated by $g_i$'s is of index $1$ or $p$ in $G_1(\infty)$ and by Sylow theorems 
normal. If the index is $p$ then the quotient is abelian, therefore $[G_1(\infty),G_1(\infty)] \subset  \langle g_i \rangle$.
This implies that $Z(G_1(\infty)) \subset \langle g_i \rangle$, a contradiction \cite{MMat}.

To each $g_i$ we associate  the variables
$X_{ij}^{g_i}$, $i>j$.  Let $X$ be a function on $\mathcal{X}$ that 
reduces to $x$ on the special fibre. 
Instead of studying arbitrary liftings in $L_{s+1}(A)$ and then 
computing the effect of taking the conjugation equivalence 
on the Krull dimension, we will choose a suitable basis on $L_i$.
We choose the elements  $\{1,X^p,X^{p^2},\ldots,X^{p^{s-1}},W\}$
as a basis for the free $A$ modules $L_i$ of proposition \ref{crit-matrix-lift}.

An element  $g_i \in \mathcal{B}$ corresponds to the automorphism given by 
\begin{equation}\label{autmatlehr}
X\mapsto X + X_{21}^{g_i}, W \mapsto W + f_{g_i}(X), 
\end{equation}
where 
\[
f_{g_i}(X)=\sum_{\nu=1}^{s-1} X_{s+1,\nu}^{g_i} X^{p^\nu}.
\]
The elements $X_{\nu,\mu}^{g_i}$ for $2\leq \nu \leq s$, $1\leq \mu \leq s$ are expressed 
in terms of $X_{21}^{g_i}$ by expanding $(X+X_{21}^{g_i})^\mu$.
We also observe that the variables $X_{s+1,\nu}^{g_i}$ are also depended to $X_{21}^{g_i}$ in the following way:
Let $\Delta(f)(X,Y):=f(X+Y)-f(X)-f(Y)$. Then, according to \cite[section 5]{CL-MM} we know that 
\[
 \Delta(f)(X,Y)=R(X,Y) + (\mathrm{Id}-F) P_f(X,Y),
\]
where $P_f(X,Y)$ is uniquely characterised by 
\[
 P_f(X,Y)=(\mathrm{Id}+F+\cdots+F^n) \Delta(f) \mod (X^{{p^{s-1}}+1}),
\]
for any $n>s$. Moreover, to  every root $y$ of $Ad_f(Y)$ corresponds an automorphism 
$X \mapsto X+y$, $W \mapsto W +P_f(X,y)$. These results combined with  (\ref{autmatlehr})  express  the variables
$X_{s+1,\nu}^{g_i}$ as functions of $X_{21}^{g_i}$, since 
\[
 f_{g_i}(X)=\sum_{\nu=1}^{s} X_{s+1,\nu}^{g_i} X^{p^\nu}=P_f(X,X_{21}^{g_i}).
\]
The desired deformation ring is given as a quotient $k[X_{21}^{g_i}: 1\leq i \leq 2s]/I$, where 
$I$ is an ideal generated by all relations among the elements $X_{21}^{g_i}$.

Moreover for the the polynomial $Ad_f(Y)=\sum_{\nu=0}^{2s} a_\nu Y^{p^\nu}$ we have 
 \[Ad_f(Y)=\sum_{0 \leq \lambda \leq s} t_{i}^{p^{s-i}} Y^{p^{s-i}}+ t_i^{p^s} Y^{p^{s+i}}\] 
This means that the coefficients $a_i$ must obey the following relation: 
\begin{equation} \label{sym1}
a_{s+\lambda}=a_{s-\lambda}^{p^\lambda} 
\end{equation}
for all $\lambda=1,\ldots,s$. Since $X_{21}^{g_i}$ are all roots of $Ad_f(Y)$ we see that the  condition  (\ref{sym1}) on coefficients
implies conditions on $X_{21}^{g_i}$. In order to describe them we need to recall some theory on additive polynomials.
Set 
\[
\Delta(x_1,\ldots,x_r)=\det 
\begin{pmatrix}
x_1 &  \cdots  & x_r \\
x_1^p & \cdots & x_r^p \\
\vdots & &  \vdots \\
x_1^{p^{r-1}} & & x_r^{p^{r-1}}
\end{pmatrix}.
\] 
We will need the following 
\begin{pro}
A  polynomial  is additive if and only if the set of its roots is  an $\mathbb{F}_p$-vector space. 
If $x_1,\ldots,x_r$ form a basis of an $\mathbb{F}_p$-vector space then a monic  additive polynomial with 
these roots is of the form:
\[
 f(Y)= \frac{\Delta(x_1,\ldots,x_r,Y)}{\Delta(x_1,\ldots,x_r)}.
\]
\end{pro}
\begin{proof}
 \cite[th. 1.2.1,lemma 1.3.6]{GossBook}
\end{proof}
The additive polynomial $Ad_f(Y)$ is given by  
\[
Ad_f(Y)=c\cdot\frac{\Delta(X_{21}^{g_1},\ldots,X_{21}^{g_{2s}},Y)}{\Delta(X_{21}^{g_1},\ldots,X_{21}^{g_{2s}})}, 
\]
where $c\neq 0$.
For the coefficients $a_i$ of $Ad_f(Y)$ we have the following formulas:
\[
 a_{2s}=c \cdot \frac{\Delta(X_{21}^{g_1},\ldots,X_{21}^{g_{2s}})}{\Delta(X_{21}^{g_1},\ldots,X_{21}^{g_{2s}})}=c.
\]
\[
 a_0=c\cdot \frac{\Delta(X_{21}^{g_1},\ldots,X_{21}^{g_{2s}})^p}{\Delta(X_{21}^{g_1},\ldots,X_{21}^{g_{2s}})}=
c\cdot \Delta(X_{21}^{g_1},\ldots,X_{21}^{g_{2s}})^{p-1}.
\]
The condition $a_{2s}=a_0^{p^s}$, implies that 
\begin{equation} \label{coef0}
c=c^{p^s} \Delta(X_{21}^{g_1},\ldots,X_{21}^{g_{2s}})^{p^s(p-1)}.
\end{equation}
For $0 < \lambda < s$, equation $ a_{s+\lambda}=a_{s-\lambda}^{p^\lambda}$ implies
\begin{equation} \label{coefl}
c\cdot \det
\begin{pmatrix}
 X_{21}^{g_1}  & \cdots & X_{21}^{g_{2s}} \\
 \vdots & & \vdots \\
(X_{21}^{g_1})^{p^s} & \cdots & (X_{21}^{g_{2s}})^{p^s} \\
\vdots & & \vdots \\
 (X_{21}^{g_1})^{p^{s+\lambda-1}} & \cdots & (X_{21}^{g_{2s}})^{p^{s+\lambda-1}}\\
(X_{21}^{g_1})^{p^{s+\lambda+1}} & \cdots & (X_{21}^{g_{2s}})^{p^{s+\lambda-1}}\\
\vdots & &\vdots\\
(X_{21}^{g_1})^{p^{2s}} & \cdots & (X_{21}^{g_{2s}})^{p^{2s}}
\end{pmatrix}=
c^p\cdot \det
\begin{pmatrix}
 (X_{21}^{g_1})^{p^\lambda}  & \cdots & (X_{21}^{g_{2s}})^{p^\lambda} \\
\vdots & &\vdots\\
(X_{21}^{g_1})^{p^{s-1}} & \cdots & (X_{21}^{g_{2s}})^{p^{s-1}}\\
(X_{21}^{g_1})^{p^{s+1}} & \cdots & (X_{21}^{g_{2s}})^{p^{s+1}}\\
 \vdots & & \vdots \\
(X_{21}^{g_1})^{p^{s+\lambda}} & \cdots & (X_{21}^{g_{2s}})^{p^{s+\lambda}} \\
\vdots & & \vdots \\
 (X_{21}^{g_1})^{p^{2s+\lambda}} & \cdots & (X_{21}^{g_{2s}})^{p^{2s+\lambda}}
\end{pmatrix}.
\end{equation}
We are now going to look at  the problem of finding the Krull  dimension  from the infinitesimal point of view.
Write $X_{21}^{g_i} \mod \epsilon^2=x_i+\epsilon A_i$, where $x_i$ is a basis for the 
vector space of the roots of $Ad_f(Y) \mod \epsilon=Y^{p^{2s}}+Y$. According to \ref{cor123}
we have to compute the dimension of the vector space in the coordinates  $A_i$.
Observe that since $p>2$  for every $a,b\in k$ we have $(a+\epsilon b)^p=a^p+\epsilon ^p b^p=a^p \mod \epsilon^2$.
Therefore (\ref{coef0}) does not give us any information, while (\ref{coefl})  is transformed to 
\[
\det
 \begin{pmatrix}
  x_1+\epsilon A_1  & \cdots & x_{2s}+\epsilon A_{2s} \\
 \vdots & & \vdots \\
x_1^{p^s} & \cdots &  x_{2s}^{p^s} \\
\vdots & & \vdots \\
 x_1^{p^{s+\lambda-1}} & \cdots & x_{2s}^{p^{s+\lambda-1}}\\
x_1^{p^{s+\lambda+1}} & \cdots & x_{2s}^{p^{s+\lambda-1}}\\
\vdots & &\vdots\\
x_1^{p^{2s}} & \cdots & x_{2s}^{p^{2s}}
 \end{pmatrix}=0 \Rightarrow 
\det
\begin{pmatrix}
   A_1  & \cdots  &A_{2s} \\
 \vdots & & \vdots \\
x_1^{p^s} & \cdots &  x_{2s}^{p^s} \\
\vdots & & \vdots \\
 x_1^{p^{s+\lambda-1}} & \cdots & x_{2s}^{p^{s+\lambda-1}}\\
x_1^{p^{s+\lambda+1}} & \cdots & x_{2s}^{p^{s+\lambda-1}}\\
\vdots & &\vdots\\
x_1^{p^{2s}} & \cdots & x_{2s}^{p^{2s}}
 \end{pmatrix}=0 
\]
The later equation gives a linear equation in $A_1,\ldots,A_{2s}$ that has 
as solution space a $2s-1$ dimensional subspace. Let us denote this subspace 
by $V_\lambda$. Observe that all vectors  $v_\nu:=(x_1^{p^\nu},\ldots,x_{2s}^{p^\nu})$ for $\nu\neq s+\lambda$ are
in $V_\lambda$.  The vectors $(x_1^{p^\nu},\ldots,x_{2s}^{p^\nu})$ are linear independent since 
$\Delta(x_1,\ldots,x_{2s})=1\neq 0$ \cite[lemma 1.3.3 ]{GossBook}.
Therefore 
\[
 V_\lambda=\langle v_{\nu}: 1 \leq \nu \leq 2s, \nu\neq \lambda+s  \rangle.
\]
All conditions we obtain from all comparisons of coefficients form the vector space
\[
 V:=\cap_{1\leq \lambda\leq s}  V_\lambda,
\]
 that is clearly a vector space of dimension $s$. This is the desired Krull dimension.

\providecommand{\bysame}{\leavevmode\hbox to3em{\hrulefill}\thinspace}
\providecommand{\MR}{\relax\ifhmode\unskip\space\fi MR }
\providecommand{\MRhref}[2]{%
  \href{http://www.ams.org/mathscinet-getitem?mr=#1}{#2}
}
\providecommand{\href}[2]{#2}

\end{document}